\documentclass[a4paper]{article}
\usepackage{a4wide}
\usepackage{authblk}
\usepackage{setspace}
\usepackage{endfloat}

\usepackage{mathrsfs,amsmath,amssymb,ascmac,cases,bm}
\usepackage{mathtools}
\usepackage{braket}

\newcommand{\nn}[0]{\notag\\}
\newcommand{\pd}[2]{\frac{\partial #1}{\partial #2}}
\newcommand{\od}[2]{\frac{d #1}{d #2}}
\newcommand{\spd}[2]{\partial #1/\partial #2}
\newcommand{\sod}[2]{d #1/d #2}
\newcommand{\inv}[1]{{#1}^{-1}}

\newcommand{\Pare}[1]{\left(#1 \right)}

\newcommand{\Angle}[1]{\left\langle#1 \right\rangle}

\newcommand{\At}[2]{\left. #1 \right\rvert_{#2}}
\newcommand{\R}[1]{\mathbb{R}^{#1}}
\newcommand{\df}[2]{\Omega^{#1}(#2)}
\newcommand{\sect}[1]{\Gamma\left({#1}\right)}
\newcommand{\bms}[1]{#1_{\#}}

\newcommand{\embr}{x}
\newcommand{\embc}{y}
\newcommand{\bi}{\vartheta}

\newcommand{\rb}[1]{E_{#1}}
\newcommand{\unitnormal}[1]{\mathcal{N}^{#1}}

\newcommand{\M}{\mathcal{M}}
\newcommand{\refer}{\mathcal{R}}

\newcommand{\metric}[1]{g[#1]}
\newcommand{\ac}[1]{\nabla[#1]}
\newcommand{\vol}[1]{\upsilon[#1]}
\newcommand{\smooth}{C^\infty}
\newcommand{\resnorm}{\mathcal{G}}
\newcommand{\bsb}[3]{B_{\Pare{{#1}, {#2}, {#3}}}} 
\newcommand{\abs}[1]{\left\lvert {#1} \right\rvert}
\newcommand{\nrasol}[0]{\tilde{\embc}}

\begin{document}

\title{Geometrical Modelling and Numerical Analysis of Dislocaion Mechanics
}
\author[1]{Shunsuke Kobayashi}
\author[1,2]{Ryuichi Tarumi}

\affil[1]{Graduate School of Engineering Science, Osaka University, 1-3 Machikaneyama-cho, Toyonaka, Osaka, 560-8531, Japan}
\affil[2]{JST PRESTO, 4-1-8 Hon-cho, Kawaguchi-shi, Saitama 332-0012, Japan}
\date{}
\maketitle
\begin{abstract}
This study undertakes the mathematical modelling and numerical analysis of dislocations within the framework of differential geometry.
The fundamental configurations, \textit{i.e.} reference, intermediate and current configurations, are expressed as the Riemann--Cartan manifold, which equips the Riemannian metric and Weitzenb\"ock connection.
The torsion 2-form on the intermediate configuration is obtained through the Hodge duality of the dislocation density and the corresponding bundle isomorphism is subjected to the Helmholtz decomposition.
This analysis introduces the boundary condition for plastic deformation.
Cartan first structure equation and stress equilibrium equation are solved numerically using weak form variational expressions and isogeometric analysis.
The numerical analysis carried out for this study reveals the distribution of plastic deformation fields around screw and edge dislocations for the first time.
It also demonstrates stress fields around dislocations of which the distant fields show full agreement with the classical Volterra theory, while at the same time eliminating the singularity otherwise introduced at the dislocation by classical methods.
The stress fields include several characteristic features
due to the geometrical nonlinearity included therein.
We also demonstrate that free surfaces affect both plastic and elastic deformation, but in different ways.
The mathematical framework of this study is applicable to an arbitrary configuration of dislocations.
\end{abstract}

\section{Introduction}

Plastic deformation of crystalline materials proceeds mainly by slip deformation between crystal planes.
Since the slip is non-uniform along the sliding direction, a linear lattice defect called a dislocation is formed on the slip plane.
Dislocation is one of the most important crystal lattice defects as it governs mechanical properties such as strength, ductility, creep resistance and fracture toughness \cite{Hirth}.
Early work on dislocations was conducted by Volterra, Orowan, Taylor and Polanyi \cite{orowan1, orowan2, orowan3,polanyi,taylor,Volterra}.
Perhaps, one of the most important early discoveries was the existence of the stress field formed around a dislocation \cite{Volterra,Mura}.
Although useful analytical formulae have been obtained, there remain several issues yet to be addressed.
One significant issue is that of the singularity contained in the current analytical expressions, in which the components of stress diverge to infinity at the core of the dislocation.
An equally important issue is the use of linear approximations.
We expect finite deformation in the vicinity of the dislocation core.
However, linear elasticity is strictly restricted to infinitesimal deformation in order to satisfy the principle of material frame indifference \cite{Marsden-Hughes}.
Modelling of dislocation based on nonlinear elasticity is essential to further progress.
Free surface is also an important issue; most of the previous studies are confined to the stress field in infinite media.
To date, many attempts to improve the modelling of dislocation stress fields have been reported.
Importantly, Lazar \textit{et al.} conducted extensive studies and succeeded removing the singularity by modifying the constitutive equation of elasticity \cite{lazar_1,lazar_2,lazar_3,lazar_4}.
Unfortunately, however, these are still far from supplying a general solution to all the current modelling problems.

Geometrical modelling of dislocations is now being pursued from a completely different perspective to that of earlier approaches.
The origin of the newer approach study dates back to the analysis of compatibility conditions by Kondo in the early 1950s.
Kondo pointed out that the degree of incompatibility corresponds to the Riemann curvature \cite{kondo_1}.
Kondo subsequently proposed a theory of dislocations based on differential geometry \cite{kondo_2}.
Independently, Bilby \textit{et al.} \cite{bilby_notitle_1955}, and Kr\"oner and Seeger \cite{kroner_nicbt-lineare_nodate} proposed an equivalent geometrical theory of dislocations.
The most important aspect of all these theories is that the incompatible state of a continuum does not exist in the ordinary Euclidean space but is treated as a mathematically generalised Riemann--Cartan manifold.
These studies also clarified the relation between the lattice defects, \textit{i.e.} dislocation and disclination, and the torsion and curvature in the affine connection \cite{kondo_2, bilby_notitle_1955, kroner_nicbt-lineare_nodate, anthony_theorie_1970}.
This theory was developed further by Noll \cite{noll}, Wang \cite{wang}, de Wit \cite{dewit}, Le and Stumpf \cite{le_stumpf_1,le_stumpf_2}, Wenzelburger \cite{wenzelburger} and Binz \textit{et al.} \cite{binz}.
Yavari and Goriely reformulated the theories using modern differential geometry and derived analytical solutions to the nonlinear boundary value problems \cite{yavari_riemanncartan_2012,yavari_riemanncartan_2013,yavari_geometry_2014}.
Similar works have been reported by Edelen \cite{edelen}, Acharya \cite{acharya} and Clayton \cite{clayton}.
However, these analyses fail to obtain the stress fields for arbitrary dislocation configurations. 
The issue arises from the mathematical method used to obtain the Riemann--Cartan manifolds.
As will be explained in later sections, it is necessary to integrate the Cartan first structure equation for a given distribution of dislocations.
Since most of the previous studies have sought analytical solutions, they have required significant restrictions upon the arrangement of dislocations in exchange.
More precisely, Willis applied the perturbation method \cite{Willis}, Yavari and Goriely used the semi-inverse method \cite{yavari_riemanncartan_2012, yavari_geometry_2014}, while Edelen \cite{edelen}, Acharya \cite{acharya} and our previous studies \cite{kobayashi_1, kobayashi_2} employed the homotopy operator.
Ortiz-Bernardin and Sfyris conducted a finite element analysis for a nonlinear elastic bar with dislocations of linear distribution \cite{ortiz-bernardin_sfyris}.
However, those researchers also employed an analytic model of the plastic deformation due to dislocations, instead of solving the Cartan first structure equation.
A highly symmetric or uniform arrangement of dislocations is still required to obtain a solution using the methods which were proposed so far.
A possible approach to a more fundamental solution is to integrate the Cartan first structure equation numerically, rather than analytically.
That method permits nonlinear analysis for an arbitrary configuration of dislocations in a finite medium and thereby expands the range of applications of the geometrical theory.

The aim of the study reported in this paper is to construct a mathematical framework for the numerical integration of the Cartan first structure equations and to implement the theory in practical terms.
The construction of this paper is as follows.
First, this introduction has furnished a brief overview of the theory of dislocations.
We have described the state-of-the-art and some issues remain to be addressed.
In the next section (Section 2), we will explain the kinematics of dislocation based on differential geometry.
The three fundamental configurations, \textit{i.e.} the reference, intermediate and current configurations, are expressed as the Riemann--Cartan manifolds.
As the primary mathematical assumption, we consider that the manifold is parallelisable.
This assumption ensures that a bundle isomorphism exists between the tangent bundle and the product bundle.
This structure enables us to properly describe the intermediate configuration as a Riemann--Cartan manifold.
We also derive the Cartan first structure equation, which governs the plastic deformation, from the Hodge duality of dislocation density and torsion of the connection.
Section 3 summarises the variational formulations, \textit{i.e.} the structure equations for plasticity and the stress equilibrium equation for elasticity.
In Section 4, we implement the weak form equations to the isogeometric analysis by numerical means, namely, a Galerkin method which uses smooth basis functions.
The results of the numerical analysis are presented in Section 5.
The analysis includes the simulation of an edge dislocation and a screw dislocation.
We compare these numerical results with the results of previous studies to verify the framework developed in this paper. 
We conclude with summary remarks in Section 6.

\section{Differential geometry of dislocation kinematics}
\subsection{Reference and current configurations}
Following the previous studies by Kondo \cite{kondo_2}, and Yavari and Goriely \cite{yavari_riemanncartan_2012}, we express the kinematics of dislocation in terms of differential geometry.
Let us first introduces the three different configurations: reference, intermediate and current.
These represent the three different states of a crystal: a stress-free perfect crystal state, a plastically deformed state due to dislocations and an elastically relaxed state following the plastic deformation.
Here the reference and current configurations are in the three-dimensional Euclidean space $\mathbb{R}^3$.
In contrast, the intermediate configuration cannot exist in $\mathbb{R}^3$.
This is because the intermediate configuration is deformed from the reference configuration in which it has a stress-free state.
On the other hand, the current configuration is obtained from the embedding of the intermediate configuration into $\mathbb{R}^3$ in such a way that the strain energy is minimised.
The embedding map corresponds to elastic deformation and its magnitude, \textit{i.e.}, elastic strain, expresses the difference in the Riemannian metric between the two configurations.
It is for this reason that differential geometry is necessary for dislocation mechanics.

Let $\M$ be a three-dimensional compact $C^\infty$ manifold with boundary.
We assume that the manifold equips a smooth Riemannian metric $g$ and $g$-compatible affine connection $\nabla$ everywhere on $\M$.
Then, the triplet $(\M, g, \nabla)$ is called a Riemann--Cartan manifold \cite{yavari_riemanncartan_2012}.
Obviously, a Riemann--Cartan manifold includes conventional Euclidean subspace as a special case possessing the Euclidean metric and flat symmetric connection.
Previous studies revealed that the Riemann--Cartan manifold is the proper mathematical framework for the dislocation mechanics \cite{yavari_riemanncartan_2012,yavari_riemanncartan_2013,yavari_geometry_2014}.

We begin with the mathematical construction of the reference configuration.
We assume that this configuration is smoothly embedded in $\mathbb{R}^3$ with the invertible $C^\infty$ map $\embr{}: \M \to \mathbb{R}^3$.
The map defines the submanifold $\M_{\mathcal{R}}=x(\M) \subset \mathbb{R}^3$ in Euclidean space.
Let $\Angle{\cdot, \cdot}_{\R{3}}\colon \R{3}\times \R{3}\to \R{}$ be the standard inner product of $\R{3}$ and let $E_1$, $E_2$, $E_3\in \R{3}$ be an orthonormal basis of $\R{3}$ such that $\Angle{E_i, E_j}_{\R{3}}=\delta_{ij}$.
Since we define the map $x$ is invertible, we have the orthonormal vector field $d\inv{x}(E_i)=\spd{}{x^i}$.
Here $d\inv{\embr}\colon \sect{T\M_\refer}\to \sect{T\M}$ is the tangent map of $\inv{\embr}$ between the set of all smooth sections $\sect{T\M_\refer}$ and $\sect{T\M}$.
Then, the map $\embr{}$ induces a Riemannian metric on $\M$ such that $\metric{\embr{}}(X, Y)=\Angle{d\embr{}(X), d\embr{}(Y)}_{\R{3}}$.
For instance, the inner product of the basis $\spd{}{x^i}$ and $\spd{}{x^j}$ is $\metric{\embr}(\spd{}{x^i}, \spd{}{x^j})=\delta_{ij}$.
Thus, we can express the Riemannian metric $\metric{\embr{}}$ as
\begin{align}\label{reference metric}
\metric{x} = \delta_{ij} dx^i\otimes dx^j,
\end{align}
where $dx^i$ is the coframe of the orthonormal basis $\spd{}{x^i}$ such that $dx^i(\spd{}{x^j})=\delta^i_j$.
A straightforward calculation yields that the Levi-Civita connection $\ac{\embr{}}$ associated with the metric $\metric{\embr{}}$ is symmetric and flat, \textit{i.e.} the connection is free from torsion and curvature \cite{tu_differential_2017}.
The triplet $(\M, \metric{\embr{}}, \ac{\embr{}})$ defines the Riemann--Cartan manifold of the reference configuration.
Hereafter, we use the frame $\partial /\partial x^i$ and coframe $dx^i$ for a local representation of tensor and differential form on the manifold $\M$.

We assume that the current configuration is also embedded in $\mathbb{R}^3$.
Let $y: \M \to \mathbb{R}^3$ be the $C^\infty$ embedding with $\M_{\mathcal{C}}=y(\M)$.
Then, we can also define the induced metric on $\M$ as $\metric{\embc{}}(X, Y) \coloneqq \Angle{d\embc{}(X), d\embc{}(Y)}_{\R{3}}$.
A local expression of $\metric{\embc{}}$ is
\begin{align}\label{eq:CurrentMetric}
  \metric{y} = \delta_{ij} \pd{y^i}{x^{k}} \pd{y^j}{x^{l}} dx^k\otimes dx^l.
\end{align}
Since the metric $\metric{\embc{}}$ is induced from the Euclidean embedding $\embc{}$, the Levi-Civita connection $\ac{\embc{}}$ is free from torsion and curvature.
The triplet $(\M, \metric{\embc{}}, \ac{\embc{}})$ defines a Riemann--Cartan manifold of the current configuration.
Note that the reference and current configuration are determined as the Riemann--Cartan manifold only by the embeddings $\embr$ and $\embc$.

\subsection{Intermediate configuration}

In contrast to the previous two cases, the mathematical construction of the intermediate configuration requires careful geometrical analysis.
Hereafter, we impose the fundamental assumption on the manifold $\M$: that it is parallelisable.
This assumption indicates that there exist smooth vector fields which can be a basis for the tangent space everywhere on $\M$ \cite{lee_introduction_2012}.
Such a vector field, $X$, is defined as a section of the tangent bundle $X\in \sect{T\M}$.
According to the standard theory of differential geometry, the tangent bundle $T\M$ and the product bundle $\M\times \R{3}$ are isomorphic over $\M$ if a base manifold $\M$ is parallelisable \cite{tu_differential_2017}.
Hence, without the loss of generality, we can introduce the trivial bundle isomorphism $\bi$ such that $\bi\colon T\M\to \M\times \R{3}$.
Mathematically, this isomorphism is understood as $\R{3}$-valued 1-form on $\M$ and, it is therefore, expressed by $\bi=\bi^i\rb{i}$ \cite{wenzelburger}.
Here, $\bi^i$ is a differential 1-form on $\M$ and $\rb{i}$ is the basis vector of $\R{3}$.
We denote the set of all $\R{3}$-valued $k$-form on $\M$ as $\Omega^k (\M; \mathbb{R}^3)$, \textit{e.g.} $\bi\in \Omega^1 (\M; \mathbb{R}^3)$.

The isomorphism $\bi$ induces the Riemannian metric $g[\bi]$ and flat affine connection $\nabla[\bi]$ on the manifold $\M$ \cite{wenzelburger}.
Let $X, Y\in \sect{T\M}$ be smooth vector fields on $\M$.
Then, the bundle isomorphism $\bi$ induces a Riemannian metric on $\M$ in such a way that $\metric{\bi}\Pare{X, Y} \coloneqq \Angle{\bms{\bi}{X}, \bms{\bi}{Y}}_{\R{3}}$.
Here $\bms{\bi}\colon \Gamma(T\M)\to \Gamma(\M\times \R{3})$ is the smooth map associated with the bundle isomorphism, $\bi$, and $\bms{\bi}X$ is the smooth section of $\M\times \R{3}$ induced by $\bi$.
For instance, the transformation of a vector field $X=\spd{}{x^k}$ into $\bms{\bi}{X}=\bi^i_k \rb{i}$.
This yields a local representation of the Riemannian metric:
\begin{align}\label{eq:IntermediateMetric}
  \metric{\bi} &{}=\delta_{ij} \bi^i_k\bi^j_l dx^k\otimes dx^l.
\end{align}
Similarly, we can also induce the affine connection using the bundle isomorphism.
It is known that a product bundle equips a trivial connection $D: \Gamma(T\M)\times \Gamma(\M\times\R{3})\to \Gamma(\M\times\R{3})$ \cite{wenzelburger}.
Let $X \in \Gamma (T\M)$ be a smooth vector field on $\M$ and let $s=s^i\rb{i}\in \Gamma(\M\times\R{3})$ be a smooth section on the product bundle $\M \times \mathbb{R}^3$.
A local form of the trivial connection is then expressed by $D_X s \coloneqq ds^i (X) \rb{i} = X^j \pd{s^i}{x^j}\rb{i}$.
The bundle isomorphism $\bi$ induces the affine connection $\nabla[\vartheta] : \Gamma(T\M)\times \Gamma(T\M)\to \Gamma(T\M)$ from $D$ such that $\ac{\bi}_X Y = \vartheta_\#^{-1} (D_{X} ({\bms{\bi}{Y}}))$.
The local expression is
\begin{align}
\label{connection}
  \ac{\bi}_X Y
  = X^k \bigg(Y^j \big(\inv{\bi}\big)^l_i \pd{\bi^i_j}{x^{k}}
  +\pd{Y^j}{x^{k}}\bigg)
  \pd{}{x^{j}},
\end{align}
where $(\inv{\bi})^i_j$ is the coefficient of the inverse map $\inv{\bi}$ such that $(\inv{\bi})^i_j\bi^j_k=\bi^i_j(\inv{\bi})^j_k=\delta^i_k$.
This is an affine connection $\ac{\bi}$ of the intermediate configuration.
A straightforward calculation confirms that the connection has a non-vanishing torsion and vanishing curvature \cite{wenzelburger}.
The flat, non-symmetric and $\metric{\bi}$-compatible affine connection $\ac{\bi}$ is termed the Weitzenb\"ock connection \cite{yavari_riemanncartan_2012, fernandez_weitzenbock_2011}.
The triplet $(\M, g[\bi], \nabla [\bi])$ defines the Riemann--Cartan manifold of the intermediate configuration.
It is noted that the intermediate configuration can be determined only through the bundle map.

The primary importance of the intermediate configuration is non-zero torsion in the Weitzenb\"ock connection $\nabla [\bi]$.
Generally, torsion of a connection $\nabla$ is defined by $T(X, Y)\coloneqq\nabla_X Y-\nabla_YX - [X, Y]$, where $[\cdot, \cdot]$ denotes the Lie bracket.
According to the definition and equation (\ref{connection}), we have a local form of the torsion of Weitzenb\"ock connection:
\begin{align}
  \label{eq:weitzenbock connection}
  T=&{}  \bigg(-\pd{\bi^i_k}{x^{j}}+\pd{\bi^i_j}{x^k}\bigg)
  \big(\inv{\bi}\big)^l_i \pd{}{x^{l}}\otimes dx^j\otimes dx^k.
\end{align}
This $(1,2)$-tensor field holds the anti-symmetric property; $T(X, Y)=-T(X,Y)$.
From the torsion tensor $T$, we can define a $\mathbb{R}^3$-valued 2-form on $\M$, \textit{i.e.} a map $\tau\colon \sect{T\M}\times \sect{T\M}\to \R{3}$ such that the local representation is
\begin{align}\label{torsion2form}
  \tau
     ={} \sum_{j<k} T^i_{jk}dx^{j}\wedge dx^{k}\otimes \rb{i},
\end{align}
where $T^i_{jk}=-\partial\bi^i_k/\partial x^j+\partial \bi^i_j/\partial x^k$ is the coefficient of this differential form.
Note that $\rb{i}$ in equation (\ref{torsion2form}) is the basis of $\mathbb{R}^3$.
As we see in the later sections, this representation of torsion is useful as it has the value in $\mathbb{R}^3$.

\section{Cartan first structure equation}
\subsection{Dislocation density}

As mentioned in the previous section, geometrical theory of dislocations was developed by Kondo, Bilby \textit{et al.} and Kr\"oner and Seeger \cite{kondo_2, bilby_notitle_1955, kroner_nicbt-lineare_nodate}.
The mathematical quantity which characterises the continuous distribution of dislocations is the dislocation density tensor $\alpha$ \cite{nye}.
For the sake of simplicity, we consider a smooth dislocation density $\alpha$ on the reference configuration.

Let $b=b^i \rb{i}$ be the Burgers vector field and let $n^j \delta_{jk}dx^k$ be the 1-form along the dislocation line.
Here, $n^i$ is the $i$-th component of the tangent vector of the dislocation line.
Then, the dislocation density field $\alpha \in \Omega^1 (\M; \mathbb{R}^3)$ is given by the following form \cite{nye}
\begin{align}\label{eq:dislocation_density}
\alpha= f b^i n^j \delta_{jk} dx^k\otimes \rb{i}.
\end{align}
Here, $f$ denotes a scaling coefficient.
To the best of the author's knowledge, the most important observation was made by Kondo \cite{kondo_2}.
In this seminal work, Kondo pointed out that the dislocation density $\alpha$ in solid mechanics is equivalent to the torsion 2-form $\tau$ in differential geometry in such a way that $\tau=*\alpha$ where $*$ stands for the Hodge star operator with respect to $\metric{\embr}$.
More precisely, we have
\begin{align}\label{eq:dislocation-torsion}
\tau=
*\alpha=\left(f b^in^l \epsilon_{ljk}\right) dx^j\wedge dx^k \otimes
\rb{i},
\end{align}
where $\epsilon_{ijk}$ is a permutation symbol.
Comparing the result with equation (\ref{torsion2form}) we have $T^i_{jk}=f b^in^l \epsilon_{ljk}$.
This relation indicates that the torsion 2-form $\tau$ is uniquely determined for an arbitrary distribution of dislocation densities.
The equivalence principle enables the analysis of dislocation using the framework of differential geometry.

Now we are ready to explain how to construct the intermediate configuration for a given distribution of dislocations.
As mentioned in the previous section, bundle isomorphism $\bi$ plays an important role in the construction of the intermediate configuration.
The external derivative for a $\R{3}$-valued $k$-form $\omega \in \df{k}{\M; \R{3}}$ is generally defined as a map
$d: \df{k}{\M; \mathbb{R}^3}\to \df{k+1}{\M; \mathbb{R}^3}$ such that $d (\omega^i\otimes \rb{i})=(d\omega^i)\otimes \rb{i}$.
Thus, the external derivative of the isomorphism $\bi$ yields
\begin{align}\label{exderivative}
  d\bi = d\Pare{\bi^i_j dx^{j}}\otimes \rb{i}
  = \sum_{j<k} \Pare{-\pd{\bi^i_k}{x^{j}}+\pd{\bi^i_j}{x^{k}}}dx^{j}\wedge dx^{k}\otimes \rb{i}.
\end{align}
Comparing with equations (\ref{torsion2form}) and (\ref{exderivative}), we end up with the Cartan first structure equation:
\begin{align}\label{eq:cartan}
  \tau = d\bi.
\end{align}
According to equation (\ref{eq:dislocation-torsion}), the left side of equation (\ref{eq:cartan}) is obtained from the dislocation density $\alpha$.
Hence, integration of the exterior derivative equation for a given torsion 2-form $\tau$ yields the bundle isomorphism $\bi$.
The bundle isomorphism determines the intermediate configuration, which is in turn responsible for the plastic deformation.
This explains how the intermediate configuration $(\M, g[\bi], \nabla [\bi])$ is constructed from a given distribution of dislocations.
It should be noted here that integration of $\bi$ along a closed curve $\partial S$ gives the total Burgers vector \cite{wenzelburger}
\begin{align}\label{eq:Burgers_circuit}
b[S]=\int_{S}d\vartheta^i \otimes \rb{i}=\int_{\partial S}\vartheta^i \otimes \rb{i},
\end{align}
where $S$ is a surface of $\M$ with the boundary $\partial S$.
Since $\bi$ has its value in $\R{3}$, the integral in equation (\ref{eq:Burgers_circuit}) is justified.
We also have the Bianchi identity $d\tau=0$ since $d\tau=dd\bi= 0$.
This implies that the dislocation is not terminated inside the material.

\subsection{Helmholtz decomposition and variational formulation}

The Cartan first structure equation (\ref{eq:cartan}) relates the $\mathbb{R}^3$-valued 2-form $\tau$ and the external derivative of the $\mathbb{R}^3$-valued 1-form $\bi$.
As the next step, we need to integrate the differential equation to obtain $\bi$, and therefore the intermediate configuration for a given distribution of dislocations.
In this mathematical analysis, the Helmholtz decomposition of $\bi$ plays an essential role.
Let $\omega, \eta \in \df{1}{\mathcal{R};\mathbb{R}^3}$ be $\R{3}$-valued 1-forms on $\M$.
As we summarise in Appendix \ref{inner product on r3 valued k-form}, the set $\df{1}{\M;\R{3}}$ equips the fibre metric $\Angle{\cdot, \cdot}\colon \df{k}{\M; \R{3}}\times \df{k}{\M; \R{3}}\to \df{3}{\M; \R{3}}$.
Thus, this metric introduce the $L^2$ product for $\df{1}{\M;\R{3}}$ such that
\begin{align}\label{norm k-form}
  \Pare{\omega, \eta}\coloneqq \int_\M \Angle{\omega, \eta}.
\end{align}
The set $\df{1}{\M;\mathbb{R}^3}$ can then be decomposed into a direct sum with respect to equation (\ref{norm k-form}) such that \cite{wenzelburger}
\begin{align}\label{Helmholtz}
    \df{1}{\M;\mathbb{R}^3} =dC^\infty(\M;\R{3})\oplus \mathcal{D}(\M;\R{3}),
\end{align}
where $d\smooth(\M; \R{3})$ is the set of all $\R{3}$-valued exact forms and $\mathcal{D}(\M; \R{3})$ is the set of $\R{3}$-valued dual exact forms.
These two subspaces are characterised by \cite{wenzelburger}
\begin{align}
    d\smooth(\M;\R{3})\coloneqq&{} \set{d\psi\in \Omega^1(\M;\mathbb{R}^3) \mid \psi\in \smooth(\M; \R{3})},\nn
    \label{dual exact space}
    \mathcal{D}(\M;\R{3})\coloneqq&{} \set{\omega\in \Omega^1(\M;\mathbb{R}^3) \mid \delta \omega=0, \omega\Pare{\unitnormal{}}=0}.
  \end{align}
Here $\delta = (-1)^{\dim \M(k+1)+1}*d*$ is a co-differential operator defined by $\Omega^1 (\M;\mathbb{R}^3)\to \Omega^0(\M;\mathbb{R}^3)$ and $\unitnormal{}$ is the outward-pointing unit normal vector on the boundary $\partial \M$.
This is known as the Helmholtz decomposition \cite{wenzelburger, binz, schwarz}.
Let $\psi \in C^\infty(\M; \R{3})$ be a component of the exact part and let $\Theta \in \mathcal{D}(\M;\mathbb{R}^3)$ be the other component in the non-exact part.
Then, equation (\ref{Helmholtz}) uniquely determines the decomposition such that
\begin{align}\label{eq:Cartan2}
    \bi^{} = d\psi+\Theta.
\end{align}
Note that $d\psi$ and $\Theta$ is $L^2$-orthogonal with respect to the norm (\ref{norm k-form}) in the sense that $\Pare{d\psi,\Theta}= 0$.

The Helmholtz decomposition yields two major consequences in terms of the Cartan first structure equation: irrelevance of the exact part and the boundary condition.
From the fundamental property of the subspaces (\ref{dual exact space}), the Cartan first structure equation becomes $\tau=d(d\psi+\Theta)=d\Theta$, since $d(d\psi)=0$ by definition.
It indicates that only the non-exact part $\Theta$ is responsible for the plastic deformation of dislocations.
Hence, without the loss of generality, we can set $\psi=x$ and $\bi^i=dx^i+\Theta^i$.
The mathematical requirement of the dual exact space (\ref{dual exact space}) also defines a proper boundary condition $\Theta|_{\partial\M}(\unitnormal{})=0$.
Considering these conditions, we introduce the following functional
\begin{align}\label{eq:plasticityfunctional}
    \mathcal{I}[\Theta, \lambda] = 
    \int_\M\frac12 \Angle{\tau-d\Theta, \tau-d\Theta }
    +\int_\M  \Angle{\lambda, \delta \Theta},
\end{align}
with the boundary condition $\Theta |_{\partial \M}(\unitnormal{})=0$.
Here, $\lambda\in \df{0}{\M; \R{3}}$ is the Lagrange multiplier.
In equation (\ref{eq:plasticityfunctional}), the integrand indicates the quadratic form of the residual of the Cartan first structure equation and the second is a constraint condition in order to satisfy $\delta \Theta=0$ on $\M$.
Hence, the solution $\Theta$ of the structure equation is characterised such that it minimises the functional $\mathcal{I}[\Theta, \lambda]$.

\subsection{Local coordinate representation}

In the following analysis, a local coordinate representation is provided for the structure equation expressed in the variational form (\ref{eq:plasticityfunctional}).
Let us first introduce the local coordinate representations of torsion 2-form $\tau$ and the non-exact part of bundle isomorphism $\Theta$ such that
\begin{align}
  \tau ={} \sum_{j<k} T^i_{jk} dx^j\wedge dx^k\otimes \rb{i},
  \quad
  \Theta ={} \Theta^i_jdx^j \otimes \rb{i}.
\end{align}
To simplify the analysis, we express the residual 2-form of the structure equation by $c =\tau-d\Theta$.
A local form of the residual $c$ is expressed by
\begin{align}
  c ={}\sum_{j<k} \Pare{T^i_{jk}-\Pare{\pd{\Theta^i_k}{x^j}-\pd{\Theta^i_j}{x^k}}} dx^j\wedge dx^k\otimes \rb{i}.
\end{align}
The inner product (\ref{eq:plasticityfunctional}) can be expressed using $c$ such that $\Angle{c, c}=c^i\wedge *c^j \Angle{E_i, E_j}_{\R{3}}$, where $*$ is the Hodge star operator with respect to $\metric{\embr}$.
For instance, $*(c^i_{12}dx^1\wedge dx^2)=c^i_{12}dx^3$ and $*(c^i_{13}dx^1\wedge dx^3)=-c^i_{13}dx^2$, \textit{etc}.
Hence, we have
\begin{align}
  c\wedge *c =\sum_{j<k}c^i_{jk}c^l_{jk}  \vol{x}\otimes \rb{i}\otimes \rb{l},
  \quad
  \Angle{c, c}=\sum_{j<k}\delta_{il} c^i_{jk}c^l_{jk}\vol{x}.
\end{align}
Similarly, the local expression of the co-differential of the non-exact part $\delta \Theta$ is 
\begin{align}
  \delta \Theta ={}(-1)^{3(1+1)+1}*d*\Pare{\Theta^i_jdx^j}\otimes \rb{i}
    ={}-\delta^{jk }\pd{\Theta^i_j}{x^k}\rb{i}.
\end{align}
From the definition above, the condition $\delta \Theta=0$ requires that $\Theta$ be divergence-free.
Thus, the inner product of $\delta \Theta\in \df{0}{\M; \R{3}}$ and $\lambda \in \df{0}{\M; \R{3}}$ reads $\Angle{\lambda, \delta \Theta}= -\delta_{il}\delta^{jk}\lambda^l\spd{\Theta^i_j}{x^k}$.
Therefore, $\mathcal{I}$ is expressed in local coordinates as 
\begin{align}\label{eq:CartanFunctional}
  \mathcal{I}[\Theta^i_j,\lambda^i] = \int_\M \frac12\sum_{j<k}\delta_{il} c^i_{jk}c^l_{jk} \vol{\embr}
  -\int_\M \delta_{il}\lambda^l \delta^{jk}\pd{\Theta^i_j}{x^k} \vol{\embr},
\end{align}
where $\vol{\embr}=dx^1\wedge dx^2\wedge dx^3$ is the volume form of $\M$.

Finally, we demonstrate the local coordinate expression of the boundary condition $\Theta|_{\partial\M}(\unitnormal{})$.
The unit normal vector field $\unitnormal{}$ at the boundary of $\M$ can be expressed in local coordinates as $\unitnormal{}=\unitnormal{i}\spd{}{x^i}$ and the above equation can be organised as follows:
\begin{align}
  \At{\Theta}{\partial\M}\Pare{\unitnormal{}} ={} \Theta^i_j dx^j\Pare{\unitnormal{k}\pd{}{x^k}}\otimes \rb{i}
  ={} \Theta^i_j \unitnormal{j} \rb{i}.
\end{align}
Thus, for each $i$ we obtain $\Theta^i_j\unitnormal{j} = 0$ at each point in $\partial\M$.

\section{Numerical implementation}

\subsection{Isogeometric analysis}

The two equations derived in the previous sections were analysed numerically using isogeometric analysis (IGA) \cite{hughes_iga}.
Isogeometric analysis uses smooth non-uniform rational B-spline (NURBS) basis functions: this is a point of difference from the conventional finite element method.
Here, we briefly summarise the numerical implementation of isogeometric analysis.

NURBS is constructed from piecewise polynomial functions called B-spline basis functions.
Let $p$ be a polynomial degree of B-spline basis and let $\xi=(\xi_1,\xi_2,\dots,\xi_m)$ be a non-decreasing sequence of $m$ real values, which is called a knot vector.
The B-spline basis function is a set of $n=m-p-1$ piecewise polynomials $\{\bsb{i}{p}{\xi}\}_{i=1,\dots,n}$ defined on the interval $I_{\xi}=[\xi_1,\xi_m)$ \cite{les_wayne}.
The subscript $(i,p,\xi)$ denotes the element number $i$, polynomial degree $p$, and knot vector ${\xi}$, respectively.
For a given triplet $(i,p,\xi)$, the B-spline function is expressed by the Cox--de Boor recursion formula \cite{les_wayne}:
\begin{align}
\bsb{i}{p}{\xi}(t)=\frac{t-\xi_i}{\xi_{i+p}-\xi_{i}}\bsb{i}{p-1}{\xi}(t)+ \frac{\xi_{i+p+1}-t}{\xi_{i+p+1}-\xi_{i+1}}\bsb{i+1}{p-1}{\xi}(t),\quad \forall t\in[\xi_1,\xi_m).
\end{align}
Note that the $0$-th order B-spline function is $B_{(i,0,\xi)}(t)=1$ if $\xi_i \le t < \xi_{i+1}$, whereas otherwise $B_{(i,0,\xi)}(t)=0$.

NURBS basis functions are constructed from the rationalisation to $B_{(i,p,\xi)}(t)$ \cite{les_wayne}.
Let $\hat{I}=I_{\xi^1}\times I_{\xi^2}\times I_{\xi^3}$ be the unit cube and let $\{B^\alpha(t)=\bsb{i_1}{p_1}{\xi_1}(t^1)\bsb{i_2}{p_2}{\xi_2}(t^2)\bsb{i_3}{p_3}{\xi_3}(t^3) \mid t=(t^1,t^2,t^3)\in \hat{I}\}$ the three-dimensional B-spline basis functions.
If we denote the real-valued coefficients called weight as $\{w^\alpha\}$, then the NURBS basis functions $N^{\alpha}(t)$ and NURBS map $x: \hat{I} \to \R{3}$
are defined by
\begin{align}\label{eq:nurbs}
  N^{\alpha}(t) = \frac{
    w^{\alpha}B^\alpha(t)
  }{
    \sum_{\beta=1}^{n} w^{\beta} B^\beta(t)
  },
  \quad
  x^i(t)=\sum_{\alpha=1}^n N^\alpha(t) a_\alpha^{i},
\end{align}
where each $(a_\alpha^1,a_\alpha^2,a_\alpha^3)\in \R{3}$ is a control point that determines the NURBS map $x$ and then the geometry of the embedded submanifold $\M_\refer=x(\hat{I})\subset \R{3}$.
Hereafter, we use the Greek letter $\alpha, \beta, \dots$ for the index related to the NURBS basis function.

\begin{figure}[ht]
    \centering
    \includegraphics[width=\textwidth]{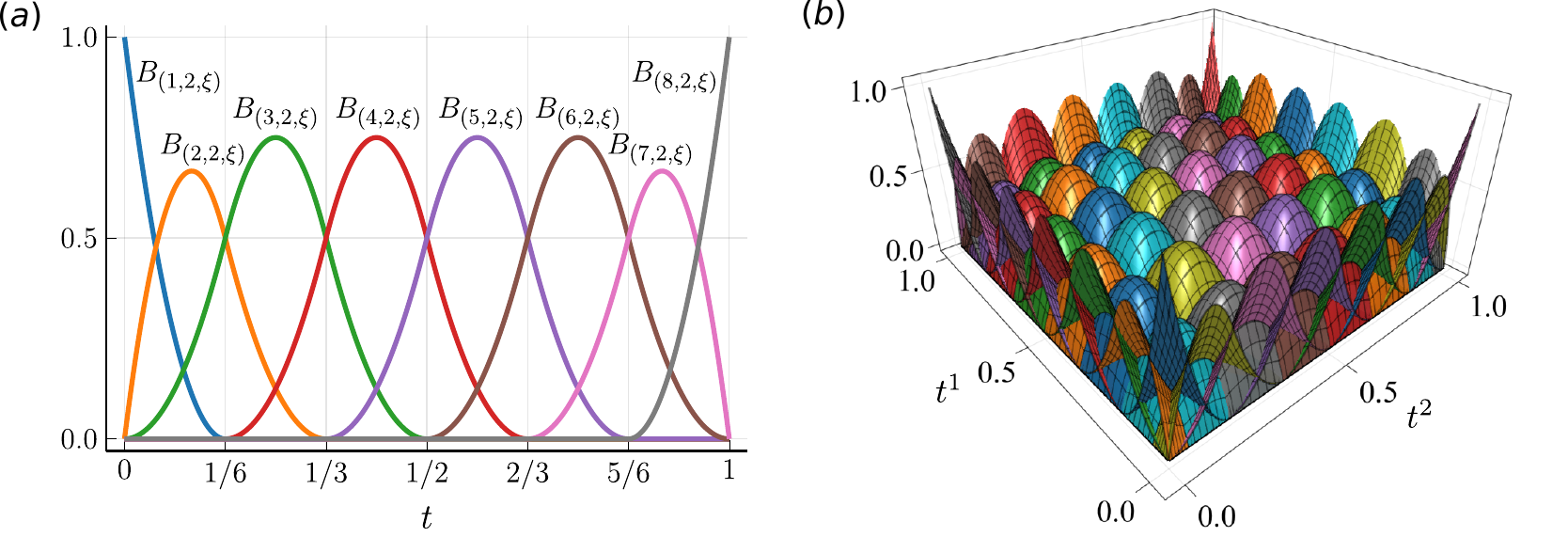}
    \caption{(a) 1-dimensional B-spline basis functions $B_{(i,2,\xi)}(t)$ of the second-order. The knot vector $\xi$ is defined by the non-decreasing sequence $\xi=(0,0,0,1/6,1/3,1/2,2/3,5/6,1,1,1)$. (b) 2-dimensional B-spline basis functions obtained by the product of $B_{(i,2,\xi)}(t)$ defined in (a). These basis functions satisfy the partitions of unity.}
    \label{fig:bsb}
\end{figure}

\subsection{Discretisation of the Cartan first structure equation}

It is clear from equation (\ref{eq:CartanFunctional}) that the functional $\mathcal{I}[\Theta^i_j,\lambda^i]$ depends on the functions $\Theta^i_j$,  and $\lambda^i$ defined on $\M$.
In this section, we apply a variational principle to $\mathcal{I}[\Theta^i_j,\lambda^i]$ and derive the discretisation of the weak form.

Now let $H$ be a matrix-valued function on $\M$ satisfying $H^i_j \unitnormal{j}=0$ and $\epsilon$ be a positive real number.
Let $c^i_{jk}(\epsilon)$ denote the function with $\Theta^i_j$ replaced by $\Theta^i_j+\epsilon H^i_j$ on the right side of equation (\ref{eq:CartanFunctional}).
Then, the first variation of equation (\ref{eq:CartanFunctional}) is obtained as follows (see Appendix \ref{Appendix:plasticity}).
\begin{align}\label{eq:variation}
  \delta \mathcal{I}[\Theta^i_j,\lambda^i] 
  =\delta_{il}\int_\M \bigg(\sum_{j<k}\od{c^i_{jk}}{\epsilon}(0) c^l_{jk}(0)
    -\eta^l \delta^{jk}\pd{\Theta^i_j}{x^k}
    -\lambda^l\delta^{jk}\pd{H^i_j}{x^k}
  \bigg) \vol{\embr}.
\end{align}
Here, the coefficient $\delta_{il}$ in the right hand side of equation (\ref{eq:variation}) suggests that we can independently determine $(\Theta^i_1, \Theta^i_2, \Theta^i_3, \lambda^i)$ for each $i=1,2,3$.
It should be noted that Acharya derives a similar weak form of the plasticity equation that corresponds to equation (\ref{eq:cartan}) \cite{acharya}.

In order to move on to the discretisaion, let us consider the $\mathbb{R}^{3\times 3}$-valued function spaces, $\mathcal{V}_m$ and $\mathcal{W}_m$, defined such that $\mathcal{V}_m =\set{\Theta^i_j \in H^1(\M)\mid \Theta^i_j \unitnormal{j}|_{\partial \M}=0}$ and $\mathcal{W}_m =\set{H^i_j \in H^1(\M)\mid H^i_j \unitnormal{j} |_{\partial \M}=0}$.
The entries $\Theta^i_j$ and $H^i_j$ are the matrix components and $\unitnormal{i}$ stands for the surface normal.
$H^1(\M)$ refers to the Sobolev space.
Similarly, we introduce the $\mathbb{R}^3$-valued function spaces, $\mathcal{V}_v$ and $\mathcal{W}_v$, defined on the reference configuration such that $\mathcal{V}_v = \set{\lambda^i\in L^2(\M)}$, and $\mathcal{W}_v ={} \set{\eta^i\in L^2(\M)}$.
The entries $\lambda^i$ and $\eta^i$ are vector components and $L^2(\M)$ is Lebesgue space.
Isogeometric analysis restricts the infinite-dimensional function spaces to finite-dimensional ones such that $\mathcal{V}_m^h$, $\mathcal{W}_m^h$, $\mathcal{V}_v^h$, and $\mathcal{W}_v^h$, respectively.
Suppose the function spaces are spanned by the NURBS basis functions.
Then, the matrix and vector components are expressed by a linear combination of $N^\alpha$ such that
\begin{align}\label{eq:nurbs coefs for plast}
  \Theta^i_j ={} \sum_{\alpha=1}^n N^\alpha (\Theta_\alpha)^i_j,\quad
  H^i_j={} \sum_{\alpha=1}^n N^\alpha (H_\alpha)^i_j,\quad
  \lambda^i={} \sum_{\alpha=1}^n N^\alpha \lambda_\alpha^i,\quad
  \eta^i={} \sum_{\alpha=1}^n N^\alpha \eta_\alpha^i.
\end{align}
Substituting the equations into the weak form of the Cartan first structure equation, we end up with the following equation:
\begin{align}\label{eq:plastic equation}
  \delta \mathcal{I}[\Theta^i_j, \lambda^i] =
  \sum_{\alpha=1}^n h_\alpha^i\left(\sum_{\beta=1}^nA^{\alpha \beta}_{ij}x_\beta^j+b^{\alpha}_i\right)
  =0.
\end{align}
Details of $A^{\alpha\beta}_{ij}, b^{i}_\alpha, h_\alpha^i, x_\beta^j$ are summarised in Appendix \ref{Appendix:plasticity}.
This equation must hold for all test functions $h^i_\alpha$.
Hence, the stationary condition $\delta \mathcal{I}[\Theta^i_j, \lambda^i]=0$ yields the equations $\sum_{\beta=1}^n A^{\alpha\beta}_{ij} x_\beta^j +b^{\alpha}_i=0$.
This is a $4n$-th order system of linear equations for the unknown coefficients $x_\beta^l$.
The solution determines the plastic deformation fields $\Theta$ due to dislocation.
Note that we use the minimum residual method (MINRES) \cite{paige} for solving equation (\ref{eq:plastic equation}) since $A^{\alpha\beta}_{ij}$ is both a symmetric and indefinite matrix.
For the convergence test of MINRES, we use the residual norm of equation (\ref{eq:plastic equation}) and we stop the MINRES iteration when the residual norm gets lower than $1.0\times 10^{-5}$.

\subsection{Stress equilibrium equation}

Finally, we present the formulation of elasticity via the geometrical theory of dislocations.
As mentioned, elastic deformation $\embc{}$ is understood as an embedding of the intermediate configuration into Euclidean space and the magnitude can be measured by the Green strain tensor $E$.
The (2,0) type tensor field can be defined as the difference of Riemann metrics between the current and intermediate configurations such that $E[\embc,\bi]= (\metric{\embc}-\metric{\bi})/2$ \cite{Marsden-Hughes,Grubic}.
According to equations (\ref{eq:CurrentMetric}) and (\ref{eq:IntermediateMetric}), a local form of the Green strain is
\begin{align}\label{eq:Green_strain}
E[\embc,\bi]= \frac{1}{2}\delta_{ij} \Pare{\pd{\embc^i}{x^k} \pd{\embc^j}{x^l}-\bi^i_k\bi^j_l}dx^k\otimes dx^l.
\end{align}
Similarly, the second Piola--Kirchhoff stress field $S$ is defined by \cite{Marsden-Hughes,Ciarlet}
\begin{align}\label{eq:2ndPK_stress}
  S[\embc, \bi] =&{} C[\bi]^{ijkl}E[\embc, \bi]_{kl} \pd{}{x^i} \otimes \pd{}{x^j},
\end{align}
where $C[\bi]^{ijkl}$ denotes the elastic coefficients.
In this study, we employ the St.Venant--Kirchhoff type constitutive equation.
The hyperelastic material is isotropic and satisfies the material frame indifference, or objectivity, for geometrically nonlinear deformation.
The elastic coefficients are given by
\begin{equation}
    C[\bi]^{ijkl}=\Pare{\frac{2\mu \nu}{1-2\nu}} \metric{\bi}^{ij}\metric{\bi}^{kl}+\mu\big(\metric{\bi}^{ik}\metric{\bi}^{jl}+\metric{\bi}^{il}\metric{\bi}^{jk}\big),
\end{equation}
where $\mu$ and $\nu$ are the shear modulus and Poisson ratio.
The metric $\metric{\bi}^{ij}$ is the inverse matrix of the coefficient $\metric{\bi}_{ij}$ in equation (\ref{eq:IntermediateMetric}).
The strain energy density is defined by $\mathcal{W}[\embc, \bi] = C[\bi](E[\embc, \bi], E[\embc, \bi])/2$.
Integration of the strain energy density over $\M$ defines the strain energy functional such that
\begin{align}\label{eq:StrainEnergyFunctional}
  W[\embc,\bi] = \int_\M \frac{1}{2} 
  C[\bi]^{ijkl} E[\embc,\bi]_{ij} E[\embc,\bi]_{kl}\upsilon[\bi],
\end{align}
where $\upsilon[\bi]=(\det\bi) dx^1\wedge dx^2\wedge dx^3$ is called the volume form.
Since the map $x\colon \M\to \R{3}$ is smoothly invertible, $W[\embc,\bi]$, the integration can be conducted in the reference configuration $\M_\refer\subset \mathbb{R}^3$.
More precisely, the functional $\delta W[\embc{}]$ depends on the bundle isomorphism $\bi$, in addition to the elastic embedding $\embc{}$.
However, as explained above, the function $\bi$ is obtained by solving the Cartan first structure equation of equation (\ref{eq:plastic equation}) numerically.
Hence, for a given intermediate configuration $\bi$, we determine the current configuration $\M_\mathcal{C}=\embc{}(\M)$ such that the embedding map $\embc{}$ minimises the strain energy functional (\ref{eq:StrainEnergyFunctional}).
This indicates that we can derive the governing equation, \textit{i.e.} the stress equilibrium equation, from the variational principle of elasticity.
We observe that the inhomogeneous distribution of plastic deformation $\bi$ is responsible for the elastic deformation $\embc{}$.
Therefore, elastic deformation occurs even if no external force is applied.

The boundary of the reference configuration $\partial \M$ can be divided into two disjoint sub-regions: the Dirichlet boundary $\Gamma_D$ and Neumann boundary $\Gamma_N$.
Let $\mathcal{V} =\set{\embc{}^i\in H^1(\M)\mid \embc{}^i |_{\Gamma_D}=D^i}$ be the function space of all admissible functions on $\M$ and let $\mathcal{W} ={} \set{h^i\in H^1(\M)\mid h^i|_{\Gamma_D}=0}$ be that of the test function.
Here, $H^1(\M)$ is the Sobolev space, $D^i$ is the prescribed boundary value on the current configuration $\M$, and $h$ is the test function satisfying the boundary condition $h^i|_{\Gamma_D}=0$.
Then, the first variation of the strain energy functional $W[\embc]$ becomes
\begin{equation}\label{eq:WFequilibriumequation}
  \delta W[\embc{}] = \int_\M
  C[\bi]^{ijkl} \delta_{mn}\pd{h^m}{x^i}\pd{\embc{}^n}{x^j} E[\embc{}, \bi]_{kl} \upsilon[\bi].
\end{equation}
This is called the stress equilibrium equation in a weak form.
This equation governs the elastic deformation in the geometrical theory of dislocations.
We solve the weak form equation numerically using isogeometric analysis.
Let us first approximate the unknown elastic deformation $\embc{}^i$ and test function $h^i$ by a linear combination of the NURBS basis functions such that 
\begin{align}
\embc{}^i=\sum_{\alpha=1}^n N^\alpha \embc{}_\alpha^i,
\quad
h^i=\sum_{\alpha=1}^n N^\alpha h_\alpha^i.
\end{align}
Inserting the approximate solutions into equation (\ref{eq:WFequilibriumequation}), we end up with the following form:
\begin{align}\label{IGA equilibrium}
  \delta_{mn}
  \int_{\hat{I}}
  C[\bi]^{ijkl}
  \pd{N^\alpha}{x^i}
  \pd{N^\beta}{x^j}
  \embc{}^n_\beta
  \frac{\delta_{pq}}{2}
  \bigg(
    \pd{N^\gamma}{x^k}
    \pd{N^\delta}{x^l}
    \embc{}^p_\gamma
    \embc{}^q_\delta
    -
    \bi^p_i
    \bi^q_j
  \bigg)
  \det{\bi}
  \det{J}
  \hat{\upsilon}=0,
\end{align}
where $\det{J}$ is the Jacobian determinant for the NURBS map $\hat{I}\to \M_\refer$ and $\hat{\upsilon}=dt^1\wedge dt^2\wedge dt^3$ is the volume form of the parameter space $\hat{I}$.
This is the system of nonlinear algebraic equations for $3n$ coefficients $\set{\embc{}^i_\alpha}_{\alpha=1,\dots,n, i=1,\dots,3}$.
To solve the nonlinear equations, we apply the Newton--Raphson method.
As explained in Appendix \ref{newton raphson method for stress equilibrium equations}, we can obtain the numerical solution of equation (\ref{IGA equilibrium}) by iteratively solving the linearised equations (\ref{eq:increment eq for see}).
When solving the linear equations, we employ the preconditioned conjugate gradient (PCG) method with the two-level overlapped additive Schwarz preconditioner for IGA \cite{BDV_iga}.
Thereby, we can reduce the number of iterations until the numerical solution converges within an acceptable tolerance.
In this paper, we stop the iterations of the Newton--Raphson and PCG methods when the relative residual norm of each equation gets smaller than $1.0\times 10^{-6}$ and $1.0\times 10^{-5}$, respectively.

\section{Results and discussion}
\subsection{Plastic deformation fields around dislocations}

\begin{figure}[htbp]
  \centering
  \includegraphics[width=\textwidth]{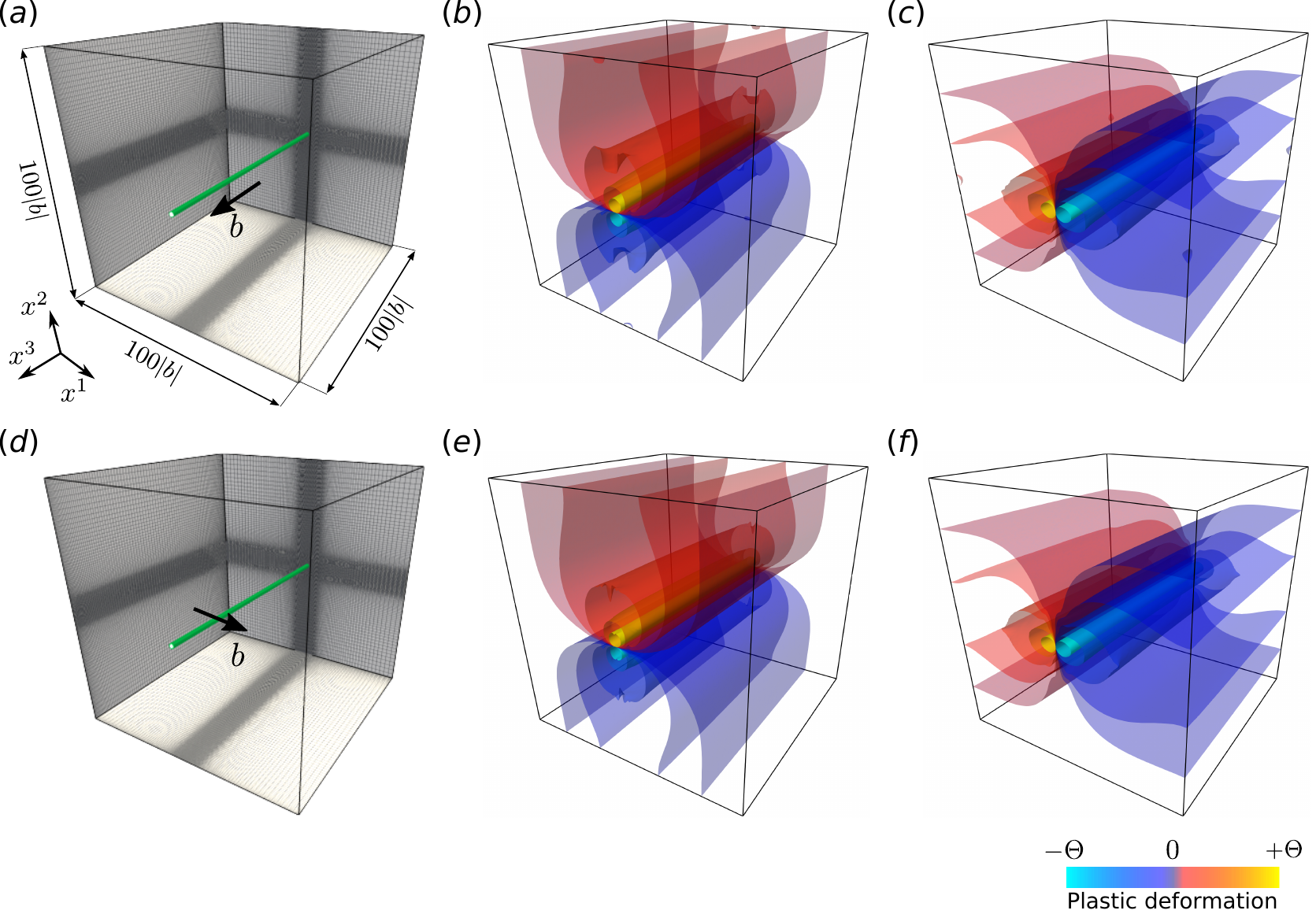}
  \caption{
(a) Schematic illustration of the screw dislocation model.
A straight dislocation line is placed at the centre of the cubic model.
Plastic deformation fields, obtained by numerically solving the Cartan first structure equation, are plotted in (b) $\Theta^3_1$ and (c) $\Theta^3_2$.
Similarly, (d) the edge dislocation model and resulting plastic deformation fields of (e) $\Theta^1_1$ and (f) $\Theta^1_2$.
The colour legend is normalised to $\Theta/b=2.9$.
Each contour plot includes 12 contour planes with the equal division to $\Theta/b$.
The results show that the $x^1$ and $x^2$ surfaces affect the distribution of plastic deformation $\Theta^i_j$ significantly, whereas little change is observed near the $x^3$ surfaces.
  }
  \label{fp isosurface}
\end{figure}

The main feature of the present study is the combined use of Helmholtz decomposition and numerical calculations for the analysis of the Cartan first structure equation.
This allows mechanical analysis for arbitrary dislocation configurations and introduces the boundary effect in the plastic deformation field.
In this section, numerical analyses are conducted for the screw and edge dislocations.
Figure \ref{fp isosurface}(a) illustrates the model and configuration of a screw dislocation.
The model has a cubic shape in the reference configuration with the dimensions of $100 |b|\times 100 |b| \times 100 |b|$.
The dislocation line is placed at the centre of the $x^1$--$x^2$ plane and parallel to the $x^3$-axis.
We use $n=200\times 200\times 100$ NURBS basis functions with uniform weight $w^\alpha=1$ to represent the cubic model.
We employ the second-order B-spline bases with non-uniform knots: \textit{i.e.} since the stress fields of dislocations tend to localise around the dislocation line, we concentrate the density of the basis functions around the dislocation line.

As explained in equation (\ref{eq:dislocation_density}), a local form of the dislocation density is given by $\alpha= f b^i n^j \delta_{jk}dx^k \otimes \rb{i}$.
According to figure \ref{fp isosurface}(a), we have $(b^1,b^2,b^3)=(0,0,b)$ and $(n^1,n^2,n^3)=(0,0,1)$ for the screw dislocation.
Similarly, we have $(b^1,b^2,b^3)=(b,0,0)$ and $(n^1,n^2,n^3)=(0,0,1)$ for the edge dislocation (see figure \ref{fp isosurface}(d)).
Consequently, the dislocation density can be written by the $\mathbb{R}^3$-valued 1-forms such that
\begin{align}\label{eq: dislocation densities}
    \alpha^\mathrm{screw}=f b dx^3 \otimes \rb{3},
    \quad
    \alpha^\mathrm{edge}=f b dx^3 \otimes \rb{1}.
\end{align}
Here, $f$ defines the radial distribution of dislocation density in a $x^1$--$x^2$ plane.
Let $r$ be a distance from the dislocation centre, \textit{i.e.} the centre of $x^1$--$x^2$ plane, and let $R$ be a radius of the dislocation core.
In this study, we employed the following radial distribution function
\begin{align}\label{eq:radial_distribution}
    f(r)=\left\{
        \begin{array}{ll}
            \dfrac{3}{\pi R^2}\left( 1-\dfrac{r}{R} \right)&\quad(r \le R)\\
            0&\quad(r>R)
        \end{array}
    \right..
\end{align}
The coefficient $\frac{3}{\pi R^2}$ is determined so that the integration of $f$ on a disk with the radius $r\geq R$ is normalised to be unity.
For an edge dislocation, therefore, the integration (\ref{eq:Burgers_circuit}) on the disk result in $b[S]^1=b^1$ while $b[S]^1< b^1$ on a disk with radius $r<R$.
The similar discussion holds for a screw dislocation.
The function $f(r)$ indicates that dislocation density distributes continuously within $r\leq R$ while becoming zero at $r=R$ and vanishing at $r>R$.
Hereafter, we call the domain $r\leq R$ the dislocation core.
The distribution of dislocation density is uniform for the dislocation line direction $x^3$.
$R$ is a characteristic length scale that determines the size of the dislocation core of the present dislocation model and the limit $R \to 0$ approaches the classical Volterra dislocation.
To simplify the analysis, we set $R=b=1$ to investigate the core of dislocation.
The relatively small size of the dislocation core is crucial when comparing and validating the results with the classical Volterra dislocations.

We obtained a numerical solution of the Cartan first structure equation for the given distribution of dislocations.
Figures \ref{fp isosurface}(b) and (c) summarise the results obtained for the screw dislocation.
The non-zero plastic deformation fields predicted from the dislocation density used in the analysis are $\Theta^3_1$, $\Theta^3_2$, and $\Theta^3_3$.
However, only $\Theta^3_1$ and $\Theta^3_2$ are plotted, as $\Theta^3_3$ was smaller than the other two components by two orders in magnitude.
The results of this analysis clearly show that the plastic deformation fields are concentrated in the vicinity of the dislocation line.
However, although the dislocation density is non-zero inside the dislocation core, they do not exist only within the core radius $R$; they smoothly distribute outside the dislocation core and some of them are observed to reach the boundary.
In addition, in the vicinity of the surface, the isosurfaces of the plastic deformation fields can be seen to be distributed so that they are perpendicular to the surface.
Essentially identical results have been obtained for the edge dislocation (see figures \ref{fp isosurface}(e) and (f)).

\begin{figure}[htbp]
  \centering
  \includegraphics[width=\textwidth]{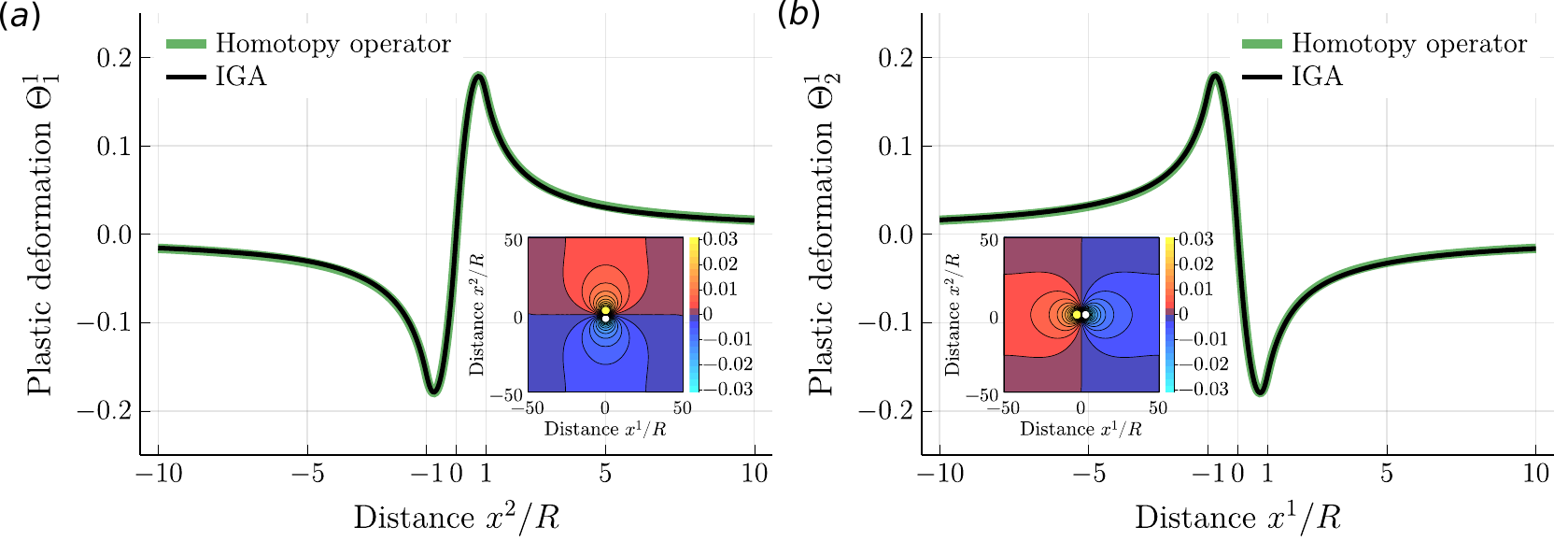}
  \caption{
    Distribution of plastic deformation obtained from the edge dislocation model.
    (a) $\Theta^1_1$ and (b) $\Theta^1_2$ show the distribution of a cross-section at the $x^3=0$ plane, \textit{i.e.} the centre of the cubic model.
    Because of the plane strain condition and the bilateral symmetry of the dislocation density, the free surfaces at $x^2/R=\pm 50$ have no effects on the plastic deformation; it enables the validation using the homotopy operator method.
    The present numerical analysis demonstrates complete agreement with analytical solutions obtained by the homotopy operator method.
  }
  \label{fp line section}
\end{figure}

To check our numerical results, we make a comparison with the results obtained, analytically, by the homotopy operator.
As has been explained, however, the homotopy operators can not consider the free boundary condition.
Hence, the comparison requires some ingenuity.
We considered the cross-section of the centre of the cube, \textit{i.e.} the $x^1$--$x^2$ plane at $x^3=0$.
Given the symmetry of the dislocation lines and dislocation density used in this study, this plane can be considered a two-dimensional plastic strain state.
In addition, as seen in figure \ref{fp isosurface}(e), the plastic deformation $\Theta^1_1$ is symmetrical on the $x^1=0$ line in this plane.
It is therefore considered that the distribution of $\Theta^1_1$ on the line is not affected by the boundary condition.
In other words, on this straight line, the analytical solution using the homotopy operator and the numerical solution using the Helmholtz decomposition are considered equivalent.
The results of this quantitative comparison for the edge dislocation are summarised in figure \ref{fp line section}.
It is clear from this result that the distributions of (a) $\Theta^1_1$ and (b) $\Theta^1_2$ obtained from the two analyses are in perfect agreement.
As shown in figure \ref{fp line section}, identical results are obtained for the other component $\Theta^1_2$ of the edge dislocation and the two components, $\Theta^3_1$ and $\Theta^3_2$, of the screw dislocation.

It is known that the elastic deformation of continua is greatly influenced by the free surface, but little analysis has been carried out on the influence of a free surface on plastic deformation.
Some previous studies have analysed the plastic deformation field by using analytical methods \cite{acharya, yavari_riemanncartan_2012}.
However, the influence of boundaries on the plastic deformation fields has not yet been considered.
In response to this problem, the present study successfully analyses the plastic deformation field in the interior of a continuum with free boundaries for the first time, by combining Helmholtz decomposition and numerical calculations.
To the best of the authors' knowledge, this is the first visualisation of the plastic deformation field formed around a dislocation with a free boundary.

\subsection{Stress fields around screw dislocations}

It is straightforward to obtain the elastic stress fields around a dislocation from the weak form stress equilibrium equation given in equation (\ref{eq:WFequilibriumequation}).
The first example of such elasticity analysis undertaken here was for the screw dislocation.
We imposed traction-free Neumann boundary conditions on the six cubic planes.
To suppress the translation and rigid body rotation, displacement of some of the cubic corners is constrained, but this has a negligibly small effect on the calculation results.

\begin{figure}[htbp]
  \centering
  \includegraphics[width=\textwidth]{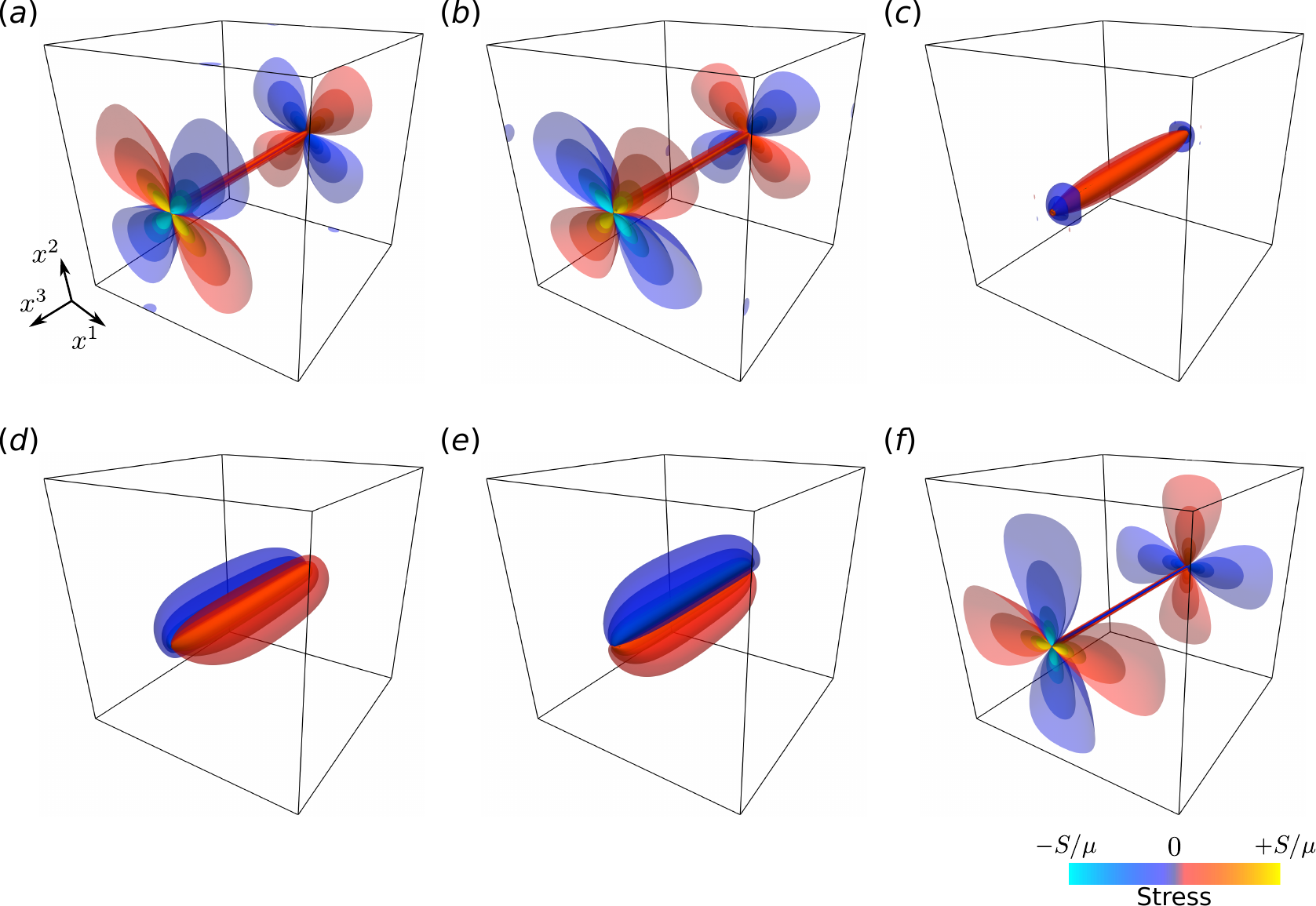}
  \caption{
   The second Piola--Kirchhoff stress fields $S^{ij}$ of the straight screw dislocation:
   (a) $S^{11}$ with $S/\mu=2.5\times 10^{-2}$,
   (b) $S^{22}$ with $S/\mu=2.5\times 10^{-2}$,
   (c) $S^{33}$ with $S/\mu=1.0\times 10^{-2}$,
   (d) $S^{23}$ with $S/\mu=5.5\times 10^{-2}$,
   (e) $S^{31}$ with $S/\mu=5.5\times 10^{-2}$ and
   (f) $S^{12}$ with $S/\mu=2.5\times 10^{-2}$.
   Each contour plot includes 12 contour planes with the equal division to $S/\mu$.
   Surface effect is responsible for the large stress concentration found in (a), (b) and (f).
   It also confirms the distribution of weak normal stress $S^{33}$ in (c) along the dislocation line.
  }
  \label{screw isosurface}
\end{figure}

Figure \ref{screw isosurface} summarises the distribution of the second Piola-Kirchhoff stress fields.
Conventional analyses, \textit{i.e.} the classical dislocation theory in an infinite medium, predicts that only $S^{23}$ and $S^{31}$ are formed around the screw dislocation.
As seen in the figure, however, the results show three distinct differences.
The first thing to state here is that the stress fields are free from singularities even at the centre of the dislocation core.
This is due to the continuous distribution of the dislocation density.
The stress fields are continuous in this analysis, because we employed the second-order NURBS basis function.
Smoothness can be improved by increasing the polynomial order of the basis function.
The second feature is the surface effect.
As predicted from the classical dislocation theory, the shear stresses $S^{23}$ and $S^{31}$ are dominant inside the material, but their magnitudes decrease as they approach the $x^3$ surfaces.
On the other hand, stress concentrations are clearly identified in $S^{11}$, $S^{22}$ and $S^{12}$ near the $x^3$ surfaces.
One possible reason for the stress concentrations is the inhomogeneity of the plastic deformation field.
As explained in the previous section, the Cartan first structure equation predicts that plastic deformation fields are disturbed by the free surface.
Another reason is compensation to satisfy the traction-free boundary condition.
These calculations imply that dislocation mechanics exhibits a size effect with the characteristic length scale of the dislocation core radius $R\,(\sim \abs{b})$.

\begin{figure}[htbp]
  \centering
  \includegraphics[width=\textwidth]{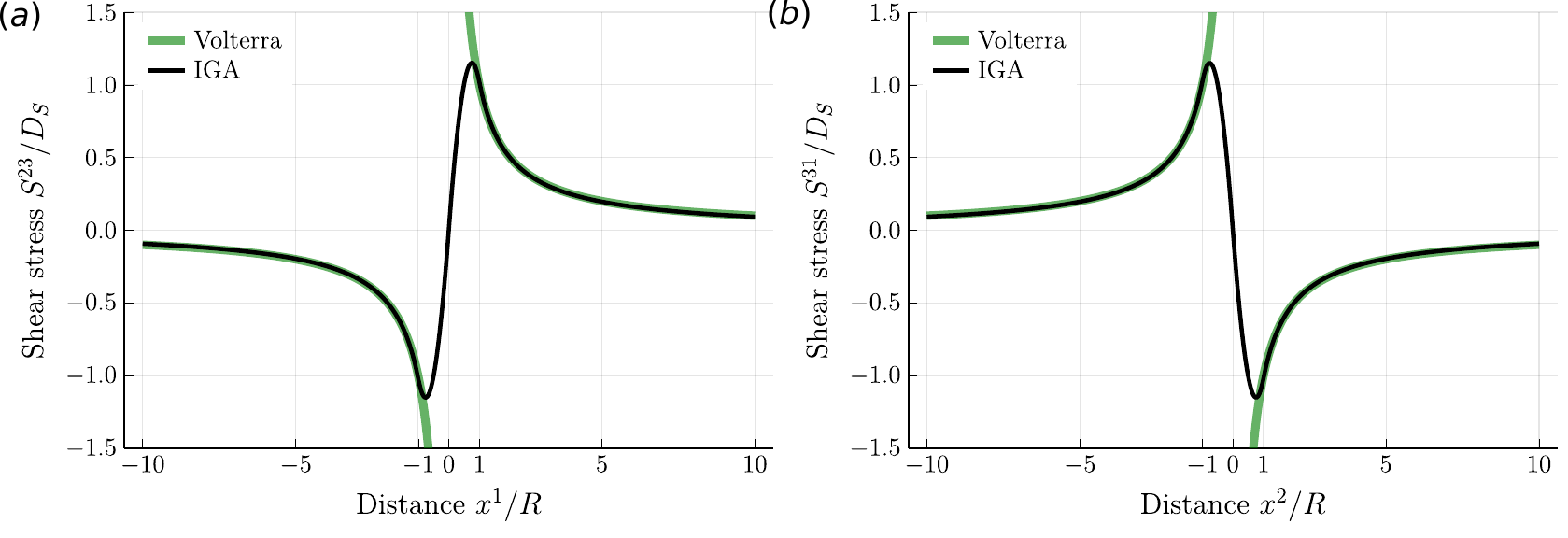}
  \caption{
    Distribution of shear stress (a) $S^{23}$ and (b) $S^{31}$ around the screw dislocation on the $x^3=0$ plane.
    Present numerical results show the complete agreement with the classical Volterra theory outside the dislocation core $x^i/R>1$.
    In addition, the present model includes no singularity even at the centre of the dislocation core.
    In this figure, the stress fields are normalised by $D_S=\mu b/2\pi R$.
  }
  \label{screw shear}
\end{figure}

Figure \ref{screw shear} shows the distribution of the dominant stress fields, $S^{23}$ and $S^{31}$, obtained at the cross-section on $x^3=0$.
As mentioned, this is a central plane of the cubic model with respect to the $x^3$ direction.
Together with the dislocation configuration given in figure \ref{fp isosurface}(a), elastic deformation in the plane is understood as the two-dimensional strain condition.
Hence, we can conduct a direct comparison to the classical Volterra theory for quantitative verification of the numerical results.
In these figures, the thin curves show the present numerical analysis while the bold ones are those obtained from the Volterra theory.
They show quantitative agreement outside the dislocation core $x^i/R>1$.
But the Volterra model diverges to infinity inside the core whereas our numerically analysed stress fields are smoothly distributed.
This result also confirmed that we were successful in removing the stress singularity of the dislocation.
Our stress fields reach their maximum around $x^i \sim R$, that is, around the edge of the dislocation core.
The stress fields are considered to approach the Volterra model in the limit $R\to 0$.
In this regard, the dislocation core radius $R$ will be a characteristic length scale of the model.

\subsection{Stress fields around edge dislocation}

\begin{figure}[htbp]
  \centering
  \includegraphics[width=\textwidth]{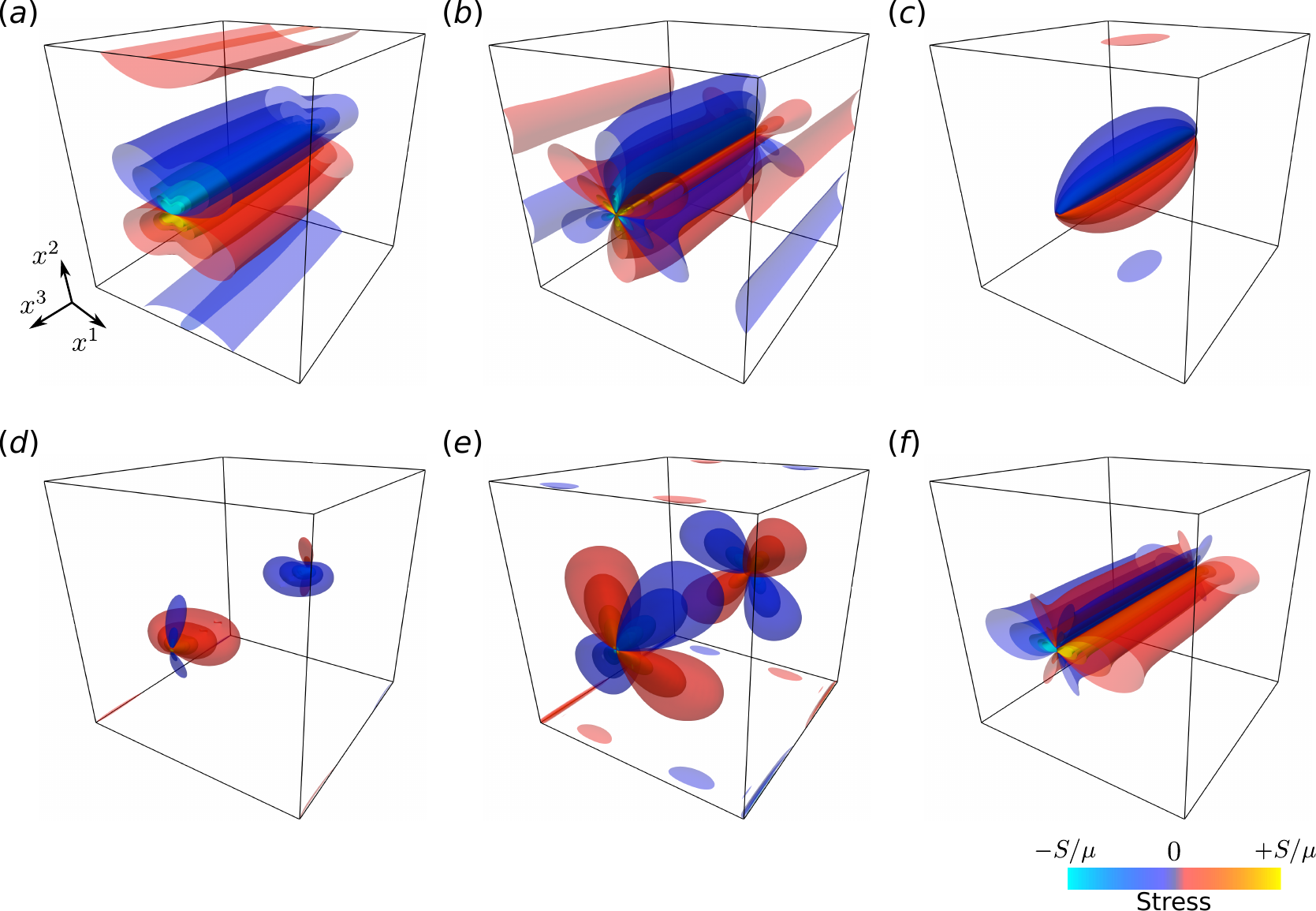}
  \caption{
   The second Piola--Kirchhoff stress fields $S$ of edge dislocation:
   (a) $S^{11}$ with $S/\mu=6.7\times 10^{-2}$,
   (b) $S^{22}$ with $S/\mu=3.8\times 10^{-2}$,
   (c) $S^{33}$ with $S/\mu=3.1\times 10^{-2}$,
   (d) $S^{23}$ with $S/\mu=1.4\times 10^{-2}$,
   (e) $S^{31}$ with $S/\mu=4.2\times 10^{-3}$ and
   (f) $S^{12}$ with $S/\mu=5.2\times 10^{-2}$.
   Each contour plot includes 12 contour planes with the equal division to $S/\mu$.
   The free surface effect is responsible for the stress concentration of $S^{23}$ and $S^{31}$ as well as the vanishing of $S^{33}$ near the $x^3$ surfaces of the cubic model.
  }
  \label{edge isosurface}
\end{figure}

We also conducted stress field analysis for the edge dislocation.
For the configuration given in figure 2(d), numerical analysis is conducted under the traction-free boundary condition.
Results of the numerical analysis are summarised in figure \ref{edge isosurface}.
As with the screw dislocation, two important properties can be identified.
At first, the stress fields shown in figure \ref{edge isosurface} are non-singular, \textit{i.e.} they do not diverge in the dislocation core.
Evidently, a continuous distribution of dislocation density $\alpha$ and plastic deformation field $\Theta^i_j$ are responsible for the non-singular stress fields.
In fact, as shown in equations (\ref{eq:Green_strain}) and (\ref{eq:2ndPK_stress}), the elastic stress field $S^{ij}$ is determined by the difference of Riemannian metric between the current and intermediate configurations.
Therefore, the non-singularity of the Riemannian metric on the intermediate configuration carries over directly to the non-singularity of the stresses.

The second characteristic feature found in figure \ref{edge isosurface} is the surface effect.
The classical dislocation theory predicts that non-vanishing stress components in the plane strain condition are $S^{11}$, $S^{22}$, $S^{33}$ and $S^{12}$ \cite{Mura}.
In view of the plane strain cross-section at $x^3=0$, the non-zero stress distributions agree with the theoretical prediction.
On the other hand, the stress distributions change significantly near the $x^3$ surfaces.
Namely, the shear stresses $S^{23}$ and $S^{31}$ appear around the dislocation core while the normal stress $S^{33}$ vanishes at the $x^3$ surfaces.
Here we take our analysis a step further by exploiting the use of differential geometry.
As mentioned above, one of the most important features of this study is the Helmholtz decomposition to the Cartan first structure equation.
This decomposition introduces the surface effects, or boundary conditions, for plastic deformation around the dislocation.
As seen in figure \ref{fp isosurface}, the surface effect on plasticity is highly anisotropic: the $x^1$ and $x^2$ surfaces affect the distribution of $\Theta^i_j$ whereas it has negligible influence on the $x^3$ plane.
In other words, in the present configuration, plastic deformation fields along the dislocation line are not affected by the boundary condition.
That approach implies that the surface stress concentration of $S^{23}$ and $S^{31}$ (or vanishing of $S^{33}$) is purely elastic, \textit{i.e.} that plastic deformation is not responsible for this phenomenon.
In this study, we employed a large-scale model because the stress field far from the dislocation core was required to validate the calculation results.
Consequently, the influence of boundary conditions on plastic deformation becomes relatively small.
However, the surface effect on plasticity would play a dominant role on the small scale such as $\sim 10R$.

\begin{figure}[htbp]
  \centering
  \includegraphics[width=\textwidth]{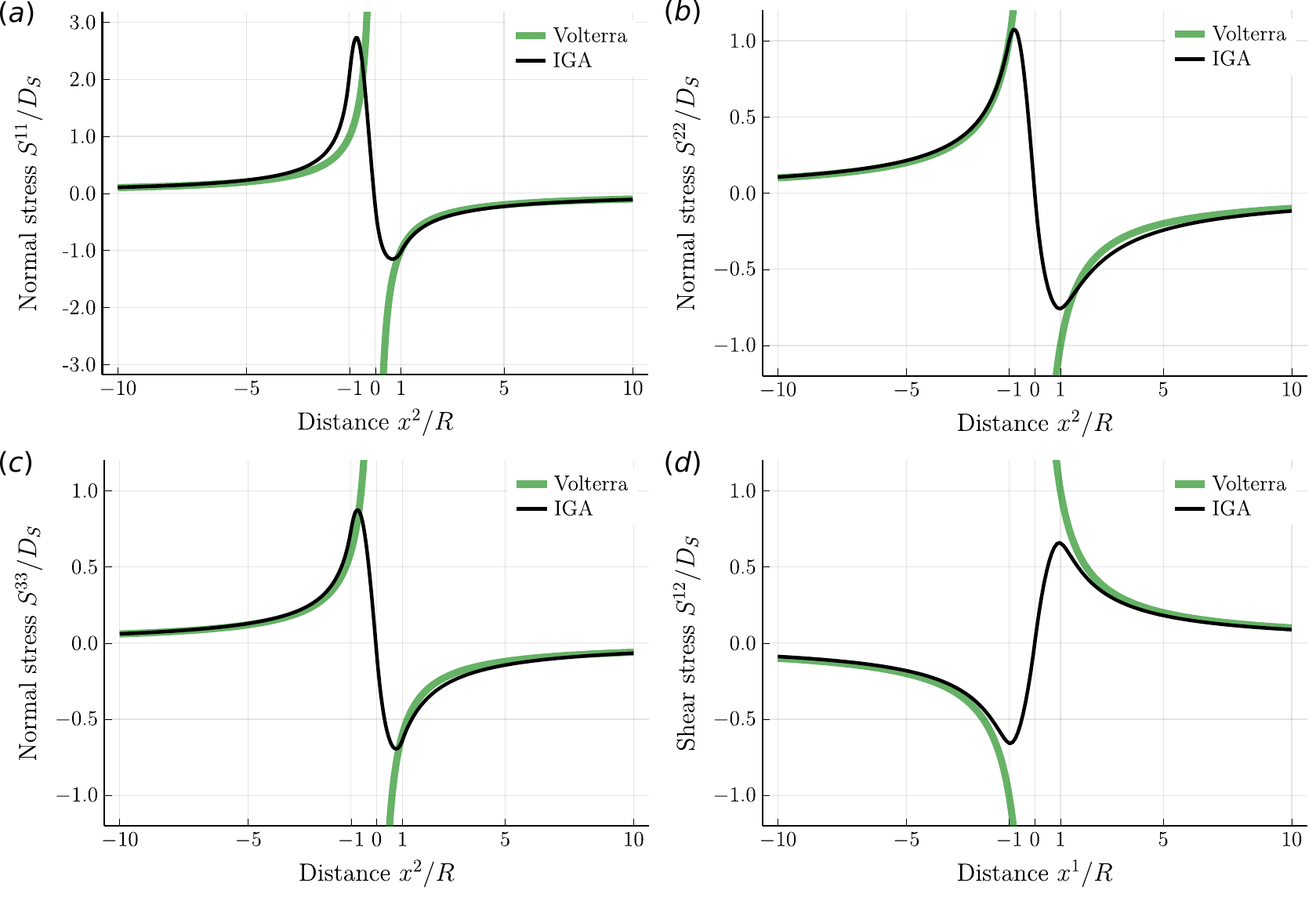}
  \caption{
  Distribution of stress (a) $S^{11}$, (b) $S^{22}$, (c) $S^{33}$ and (d) $S^{12}$ around the edge dislocation on the $x^3=0$ plane.
  The distant stress fields $x^2/R>5$ show agreement with the classical Volterra theory, whereas notable deviations are confirmed in the normal stress fields especially near the dislocation core.
  From a systematic numerical analysis with different magnitude of Burgers vector $b$, we conclude that the geometrical nonlinearity is responsible for this phenomenon.
  Note that the stress fields are normalised by $D_S=\mu b/2\pi(1-\nu) R$.
  }
  \label{edge cross section}
\end{figure}

Finally, we analysed the stress fields on the plane strain cross-section $x^3=0$ for validation of our numerical results concerning edge dislocations.
Figure \ref{edge cross section} summarises the stress fields around the core of the edge dislocation.
The stress fields far from the dislocation core, $x^2/R>5$, show agreement with the classical Volterra theory.
On the other hand, non-negligible deviations are confirmed near the dislocation core, especially in the normal stress fields.
For instance, the normal stress $S^{11}$ is asymmetric with respect to the radial distance from the dislocation core.
As seen in the analysis of screw dislocation, the geometrical nonlinearity induces normal stress fields around the dislocation core; this theory of the causation of the stress fields is supported by systematic numerical analysis.
In fact, the deviation becomes negligible when the magnitude of the Burgers vector goes to zero.
Our model thus converges to the Volterra dislocation model in the limit $b\to 0$ and $R\to 0$.

\newpage

\section{Conclusion}

Differential geometry appears to be the most appropriate mathematical framework for the theory of dislocations.
It introduces the three independent configurations, \textit{i.e.} the reference, intermediate and current configurations, as the Riemann--Cartan manifolds.
The central idea here is the equivalence of dislocation density $\alpha$ and torsion 2-form $\tau$, which leads to the celebrated Cartan first structure equation.
The main theoretical foundation of the present study was the Helmholtz decomposition of the differential equation, after the previous work by Wenzelburger.
We also conducted the numerical implementation using nonlinear elasticity and the calculus of variation.
Our numerical analysis reveals many new aspects of the mechanics of screw and edge dislocations, with superior analysis of screw and edge dislocation cores including the elimination of the singularities predicted by classical analytic methods and equivalent, validating results to the classical methods outside the dislocation cores.
The conclusions of this study are summarised as follows.
\begin{enumerate}
    \item 
    We introduced the three different configurations as Riemann--Cartan manifolds.
    The manifolds share the simply connected Euclidean subset $\M$ with a Riemannian metric $g[\cdot]$ and a Weitzenb\"ock connection $\nabla [\cdot]$.
    The intermediate configuration $(\M, \metric{\bi}, \ac{\bi})$ is characterised by the bundle isomorphism $\bi$, which is divided into two parts using the Helmholtz decomposition.
    That is, the decomposition of $\mathbb{R}^3$-valued 1-form $\Omega^1(\M;\R{3})$ into the exact form $dC^\infty(\M;\R{3})$ and dual exact form $\mathcal{D}(\M;\R{3})$.
    This is a $L^2$-orthogonal decomposition in the sense of the standard norm for $\Omega^1(\M;\R{3})$.
    The mathematical structure naturally introduces the boundary conditions for the dual exact form $\mathcal{D}(\M;\R{3})$, which turns out to be the boundary conditions for the plastic deformation due to dislocations.
    \item 
    We express the Cartan first structure equation in a variational form, \textit{i.e.} minimisation of the residual norm functional with a subsidiary condition, and we also derived the weak form equation for NURBS-based isogeometric analysis.
    The numerical implementation is advantageous as it is applicable for an arbitrary distribution of dislocations.
    Validation is conducted by cross-sectional analysis on a plane strain condition and full consistency was confirmed with results obtained using the homotopy operator.
    We were surprised to discover that plastic deformation $\Theta^i_j$ extends gradually to the free boundary even though its source, \textit{i.e.} the dislocation density $\alpha$, is strictly confined within the small dislocation core.
    This would be the first result that reveals the distribution of plastic deformation with the appropriate boundary condition.
    \item 
    Elasticity analysis was also conducted.
    Our analysis of stress fields around a screw dislocation reveals that there is no singularity even at the dislocation core.
    This is because the dislocation density is continuous in contrast to the Dirac delta function used in classical theory.
    The singularity-free stress fields show agreement with Volterra dislocation theory outside the dislocation core.
    We also reveal the emergence of normal stress fields around the dislocation core due to geometrical nonlinearity.
    On the other hand, significant localisation of normal stresses is confirmed in the vicinity of free surfaces.
    This result clearly illustrates the fact that free surfaces affect both plastic and elastic deformation, but in different ways.
    \item 
    As with the screw dislocation, our analysis of stress fields around an edge dislocation creates no singularity at the dislocation core.
    Although the distant stress fields show agreement with the Volterra dislocation theory, differences are observed in the vicinity of the dislocation core.
    In addition, the normal stress fields are non-symmetric with respect to the sign, \textit{i.e.} tension and compression are non-symmetric.
    A systematic numerical analysis revealed that geometrical nonlinearity was responsible for those non-symmetries.
\end{enumerate}

\vskip6pt

\enlargethispage{20pt}


\providecommand{\dataccess}[1]
{
  \small	
  \textbf{{Data Accessibility.}} #1
}
\providecommand{\aucontribute}[1]
{
  \small	
  \textbf{{Authors' Contributions.}} #1
}
\providecommand{\competing}[1]
{
  \small	
  \textbf{{Competing Interests.}} #1
}
\providecommand{\funding}[1]
{
  \small	
  \textbf{{Funding.}} #1
}
\providecommand{\ack}[1]
{
  \small	
  \textbf{{Acknowledgements.}} #1
}





\noindent
\ack{
  This work was supported by JST, PRESTO Grant Number JPMJPR1997 Japan and JSPS KAKENHI Grant Number JP18H05481.
  The authors gratefully acknowledge Dr A. Suzuki for discussions on numerical analysis.
  The authors would like to thank Enago (www.enago.jp) for the English language review.
  }

\appendix

\section{Mathematical foundations of the dislocation theory}
\subsection{Inner product for $\R{3}$-valued $k$-forms}\label{inner product on r3 valued k-form}
Let $\M$ be a three-dimensional smooth manifold with boundary and $g[x]$ be the induced metric defined as equation (\ref{reference metric}).
We denote the Hodge star operator with respect to $g[x]$ as $*\colon \df{k}{\M}\to \df{k-1}{\M}$.
For an orthonormal dual frame $dx^1$, $dx^2$ and $dx^3 \in \df{1}{\M}$, the Hodge star operator $*$ then acts as $*dx^1 = dx^2\wedge dx^3$, $*dx^2 = dx^3\wedge dx^1$ and $*dx^3 = dx^1\wedge dx^2$.
Similarly, the star operator $*$ acts on the $\R{3}$-valued 1-form $\omega=\omega^iE_i\in \df{1}{\M, \R{3}}$ as
\begin{align}\label{eq:hodge star for r3-valued 1-form}
  *\omega = *(\omega^i_j dx^j)E_i
  = (\omega^i_1 dx^2\wedge dx^3+\omega^i_2 dx^3\wedge dx^1+\omega^i_3 dx^1\wedge dx^2)E_i.
\end{align}
It is evident from equation (\ref{eq:hodge star for r3-valued 1-form}) that $*\omega$ is the $\R{3}$-valued 2-form on $\M$.
Now, we introduce the inner product on $\df{k}{\M;\R{3}}$ by using the star operator $*$.
Let $\Angle{\cdot, \cdot}_{\R{3}}$ be the standard inner product on the three-dimensional real vector space $\R{3}$ and let $\rb{1}$, $\rb{2}$ and $\rb{3}$ be a basis of $\R{3}$.
For any $\R{3}$-valued $k$-form $\omega=\omega^iE_i$ and $\eta=\eta^iE_i$, we can define the inner product such that \cite{wenzelburger}
\begin{align}\label{eq:inner product for df}
  \Angle{\omega, \eta} \coloneqq \omega^i\wedge *\eta^j\Angle{\rb{i}, \rb{j}}_{\R{3}}.
\end{align}
When $k=1$, for instance, we have
\begin{align}
  \Angle{\omega, \eta} = \delta_{ij} (\omega^i_1\eta^j_1 +\omega^i_2\eta^i_2 +\omega^i_3\eta^j_3)dx^1\wedge dx^2\wedge dx^3.
\end{align}
By definition, the inner product $\Angle{\cdot, \cdot}$ depends on the Riemannian metric $g[x]$ since the star operator $*$ does also.
From the definition above, this inner product can be thought of as a map $\Angle{\cdot, \cdot}\colon \df{k}{\M;\R{3}}\times \df{k}{\M;\R{3}}\to \df{3}{\M;\R{3}}$.

\subsection{First variation of the plasticity functional}\label{Appendix:plasticity}
We apply the variational principle to the functional $\mathcal{I}[\Theta,\lambda]$.
For simplicity, we denote $\resnorm[\Theta]=\int_\M \Angle{c, c}/2$ and $\mathcal{L}[\Theta, \lambda]=\int_\M \Angle{\lambda, \delta \Theta}$.
Thus, we can write $\mathcal{I}[\Theta, \lambda]=\resnorm[\Theta]+\mathcal{L}[\Theta, \lambda]$.

First, we consider the first variation of $\resnorm[\Theta]$.
Let $H$ be a matrix-valued function on $\M$ satisfying $H^i_j \unitnormal{j}=0$ and $\epsilon$ be a positive real number.
Let $c^i_{jk}(\epsilon)$ be the function with $\Theta^i_j$ replaced by $\Theta^i_j+\epsilon H^i_j$ in the right side of equation (\ref{eq:CartanFunctional}).
The first variation of the functional $\resnorm[\Theta]$ reads
\begin{align}\label{eq:variation resnorm}
  \delta \resnorm[\Theta] =&{} \At{\od{}{\epsilon}}{\epsilon=0}
    \frac12 \int_\M \sum_{j<k}\delta_{il} c^i_{jk}(\epsilon)c^l_{jk}(\epsilon) \vol{\embr}
  =
    \int_\M \sum_{j<k}\delta_{il}\od{c^i_{jk}}{\epsilon}(0) c^l_{jk}(0) \vol{\embr},
\end{align}
where $\sod{c^i_{jk}}{\epsilon}(0)$ is
\begin{align}
  \od{c^i_{jk}}{\epsilon}(0) =&{} \left.\od{}{\epsilon}\Pare{T^i_{jk}-\Pare{\pd{}{x^j}\Pare{\Theta^i_k+\epsilon H^i_k}-\pd{}{x^k}\Pare{\Theta^i_j+\epsilon H^i_j}}}\right|_{\epsilon=0}\nn
  =&{} -\Pare{\pd{H^i_k}{x^j}-\pd{H^i_j}{x^k}}.
\end{align}
Therefore, equation (\ref{eq:variation resnorm}) becomes
\begin{align}
  \label{eq:gnorm variation}
  \delta \resnorm[\Theta] =&{} \int_\M -\sum_{j<k}\delta_{il} \Pare{\pd{H^i_k}{x^j}-\pd{H^i_j}{x^k}} 
  \Pare{T^l_{jk}-\Pare{\pd{\Theta^l_k}{x^j}-\pd{\Theta^l_j}{x^k}}}
  \vol{\embr}.
\end{align}

Similarly, we show the first variation $\delta \mathcal{L}[\Theta, \lambda]$.
Let $\eta$ be any $\R{3}$-valued function on $\M$, then $\delta \mathcal{L}[\Theta, \lambda]$ can be expressed as follows:
\begin{align}
  \delta\mathcal{L}[\Theta, \lambda]=&{} \At{\od{}{\epsilon}}{\epsilon=0}\int_\M \Pare{ -\delta_{il}(\lambda^l+\epsilon \eta^l) \delta^{jk}\pd{}{x^k}\Pare{\Theta^i_j+\epsilon H^i_j}} \vol{\embr}\nn
  \label{eq:lagrange variation}
  =&{} \int_\M \Pare{
    -\delta_{il}\eta^l \delta^{jk}\pd{\Theta^i_j}{x^k}
    -\delta_{il}\lambda^l\delta^{jk}\pd{H^i_j}{x^k}
  } \vol{\embr}.
\end{align}
Therefore, we can derive the first variation of $\delta \mathcal{I}[\Theta, \lambda]$ as $\delta \resnorm[\Theta]+\delta \mathcal{L}[\Theta, \lambda]=0$.

If we take equation (\ref{eq:nurbs coefs for plast}) and substitute into equations (\ref{eq:gnorm variation}) and (\ref{eq:lagrange variation}), then we obtain the following discretisation:
\begin{align}
    \delta\resnorm[\Theta] =&{} -\delta_{il} \int_{\hat{I}} \sum_{j<k}\sum_{\alpha=1}^n \Pare{\pd{N^\alpha }{x^j}(H_\alpha)^i_k-\pd{N^\alpha}{x^k} (H_\alpha)^i_j}\nn
    &{}\times \Pare{T^l_{jk}-\sum_{\beta=1}^n \Pare{\pd{N^\beta}{x^j}(\Theta_\beta)^l_k-\pd{N^\beta }{x^k}(\Theta_\beta)^l_j}}
    \det{J}
    \hat{\upsilon}, \\
    \delta\mathcal{L}[\Theta, \lambda] =&{} -\delta_{il}
    \int_{\hat{I}} \sum_{\alpha, \beta=1}^n \Biggl(
      N^\alpha\eta^l_\alpha \delta^{jk}\pd{N^\beta}{x^k} (\Theta_\beta)^i_j
      +N^\alpha\lambda^l_\alpha \delta^{jk}\pd{N^\beta}{x^k} (H_\beta)^i_j
    \Biggr)
    \det{J}
    \hat{\upsilon},
\end{align}
where $\det{J}$ is the Jacobian determinant for the NURBS map $\hat{I}\to \M_\refer$ and $\hat{\upsilon}=dt^1\wedge dt^2\wedge dt^3$ is the volume form of $\hat{I}$.
Thus, we can denote $\delta \mathcal{I}[\Theta, \lambda]=0$ as $\sum_{\alpha=1}^n h_\alpha^i(\sum_{\beta=1}^n A^{\alpha \beta}_{ij}x_\beta^j+b^{\alpha}_i)=0$, $(i,j=1,2,3,4)$.
Here, $A^{\alpha\beta}_{ij}, b^{\alpha}_i, h_\alpha^i, x_\beta^j$ are as follows:
\begin{align}
  \label{eq:optimize fp matrix}
  \begin{split}
    \left(A^{\alpha\beta}_{ij}\right)_{i,j=1,2,3,4}=&{}\int_{\hat{I}}
    \Pare{
    \begin{matrix}
      N^\alpha_{,2}N^\beta_{,2}+N^\alpha_{,3}N^\beta_{,3}& -N^\alpha_{,2}N^\beta_{,1}& -N^\alpha_{,3}N^\beta_{,1}& N^\alpha N^\beta_{,1}\\
      -N^\alpha_{,1}N^\beta_{,2}& N^\alpha_{,1}N^\beta_{,1}+N^\alpha_{,3}N^\beta_{,3} & -N^\alpha_{,3}N^\beta_{,2}& N^\alpha N^\beta_{,2}\\
      -N^\alpha_{,1}N^\beta_{,3} & -N^\alpha_{,2}N^\beta_{,3}& N^\alpha_{,1}N^\beta_{,1}+N^\alpha_{,2}N^\beta_{,2}& N^\alpha N^\beta_{,3}\\
      N^\alpha N^\beta_{,1}&N^\alpha N^\beta_{,2}&N^\alpha N^\beta_{,3}& 0
    \end{matrix}
    }\det{J}\hat{\upsilon},
    \\
    \left(b^{\alpha}_i\right)_{i=1,2,3,4}=&{}\delta_{il}\int_{\hat{I}} \Pare{
    \begin{matrix}
      T^l_{12}N^\alpha_{,1}-T^l_{13}N^\alpha_{,3}\\
      T^l_{23}N^\alpha_{,3}-T^l_{12}N^\alpha_{,1}\\
      T^l_{31}N^\alpha_{,1}-T^l_{23}N^\alpha_{,2}\\
      0
    \end{matrix},
    }\det{J}\hat{\upsilon},
    \\
    \left(h_\alpha^i\right)_{i=1,2,3,4}=&{}\Pare{\begin{matrix}
      (H_\alpha)^i_1&(H_\alpha)^i_2&(H_\alpha)^i_3& \eta_\alpha^i
    \end{matrix}
    },\\
    \left(x_\beta^j\right)_{j=1,2,3,4} =&{} \Pare{\begin{matrix}
      (\Theta_\beta)^j_1\\
      (\Theta_\beta)^j_2\\
      (\Theta_\beta)^j_3\\
      \lambda_\beta^j
    \end{matrix}
    }.
  \end{split}
\end{align}
Here, $N^\alpha_{,i}=\pd{N^\alpha}{x^i}$ denotes the partial derivative of the basis function $N^\alpha$ with respect to $i$-th coordinate.
Thus, the order of the matrix $A^{\alpha\beta}_{ij}$ is $4n$.

\subsection{Stress equilibrium equation}\label{newton raphson method for stress equilibrium equations}
In this section, we summarise the Newton--Raphson method for the system of nonlinear equations equation (\ref{IGA equilibrium}).
For simplicity, we denote the left side of equation (\ref{IGA equilibrium}) as $f[\embc{}]^{\alpha}_m$.
This notation allows us to write $\delta W[\embc{}]=f[\embc{}]^{\alpha}_mh_\alpha^m$.
Let $\nrasol=\embr$ be an initial guess at the solution.
In general, due to the plastic deformation, this initial guess will not satisfy the equation.
To obtain a better approximation for the solution of equation (\ref{IGA equilibrium}), we take an increment $\Delta \nrasol$ and determine it in the following way.
At first, we substitute $\nrasol+\Delta \nrasol$ into equation (\ref{IGA equilibrium}).
The Taylor series expansion of the equation $f[\nrasol+\Delta \nrasol]^\alpha_m=0$ around $\nrasol$ then yields the linear approximation $f[\nrasol+\Delta\nrasol]^\alpha_m
\simeq f[\nrasol]^\alpha_m+A[\nrasol]^{\alpha\beta}_{mn}\Delta \nrasol^\beta_n=0$.
Thus, the coefficients of the increment $\Delta \nrasol$ can be determined by solving the following equation:
\begin{align}\label{eq:increment eq for see}
  A[\nrasol]^{\alpha\beta}_{ij}
  (\Delta \nrasol)_\beta^j
  =-f[\nrasol]^\alpha_i.
\end{align}
Here, $A[\nrasol]^{\alpha\beta}_{mn}$ is defined by
\begin{align}
  A[\nrasol]^{\alpha \beta}_{mn}=
    \pd{f^\alpha_m}{\nrasol{}^\beta_n}[\nrasol]=&{}
    \delta_{mn}
    \int_{\hat{I}}
    C[\bi]^{ijkl}
    \Biggl(
        \delta_{pq}
        \pd{N^\alpha}{x^i}
        \pd{N^\beta}{x^j}
        \pd{N^\gamma}{x^k}
        \pd{N^\delta}{x^l}
        \tilde{\embc{}}^p_\gamma
        \tilde{\embc{}}^q_\delta
      \nn
      &{}+
      \frac{\delta_{pq}}{2}
      \bigg(
        \pd{N^\gamma}{x^i}
        \pd{N^\delta}{x^j}
        \tilde{\embc{}}^p_\gamma \tilde{\embc{}}^q_\delta
        -
        \bi^p_i
        \bi^q_j
      \bigg)
      \Pare{
        \pd{N^\alpha}{x^k}
        \pd{N^\beta}{x^l}
      }
    \Biggr)
    \det{\bi} 
    \det{J}
    \hat{\upsilon}.
\end{align}
Equation (\ref{eq:increment eq for see}) is the system of linear equation with $3n$ unknown coefficients $\set{\nrasol^i_\alpha}_{\alpha=1,\dots, n,i=1,2,3}$.
If we obtain the solution of equation (\ref{eq:increment eq for see}), the updated guess $\nrasol+\Delta\nrasol$ removes the linear terms of $f[\nrasol+\Delta\nrasol]^\alpha_m$.
This means that $\nrasol+\Delta\nrasol$ is the better guess for the solution of equation (\ref{IGA equilibrium}).
Therefore if we iterate solving equation (\ref{eq:increment eq for see}) and updating the current guess $\nrasol$, then we obtain the approximated solution for equation (\ref{IGA equilibrium}) within an acceptable tolerance.

\end{document}